\documentclass[11pt,reqno,preprint]{elsarticle}

\usepackage{placeins}

\makeatletter
\def\ps@pprintTitle{%
 \let\@oddhead\@empty
 \let\@evenhead\@empty
 \def\@oddfoot{}%
 \let\@evenfoot\@oddfoot}
\makeatother

\usepackage{amsmath,amsthm,amscd,wrapfig,amsfonts,amssymb,mathtools}
\usepackage{enumitem,comment,mathabx,mathrsfs}
\usepackage{booktabs,multirow,lscape,cool}
\usepackage{url,parskip,caption}
\usepackage[top=1.0in, bottom=1.0in, left=1.0in, right=1.0in]{geometry}

\providecommand{\e}[1]{\ensuremath{\times 10^{#1}}}

\makeatletter
\g@addto@macro\normalsize{%
  \setlength\abovedisplayskip{.4em}
  \setlength\belowdisplayskip{.4em}
  \setlength\abovedisplayshortskip{.4em}
  \setlength\belowdisplayshortskip{.4em}
}

\begin{document}

\begin{frontmatter}

\begin{abstract}

Standard solvers for the variable coefficient Helmholtz equation in two spatial
dimensions have running times
which grow quadratically with the wavenumber $k$.
Here, we describe a solver which  applies only when the 
scattering  potential is radially symmetric but whose running time is $\mathcal{O}\left(k \log(k) \right)$
in typical cases.
We also present the results of numerical experiments demonstrating the properties
of our solver,  the code for which is publicly available.

\end{abstract}


\title
{
A quasilinear complexity algorithm for the numerical simulation of scattering from
a two-dimensional radially symmetric potential
}

\author{James Bremer}
\ead{bremer@math.ucdavis.edu}
\address{Department of Mathematics, University of California, Davis}

\end{frontmatter}

\begin{section}{Introduction}

In the frequency domain, the displacement $u$ in an inhomogeneous fluid satisfies the 
variable coefficient Helmholtz equation
%
\begin{equation}
\Delta u(x) + k^2 \left( 1+q(x) \right) u(x) = 0 \ \ \mbox{for all}\  x \in \mathbb{R}^2.
\label{introduction:helm1}
\end{equation}
In many applications, the wavenumber $k$ is real-valued
and the scattering potential $q$ is piecewise smooth with compact support contained inside
of a disk $\Omega$ of radius $R$ centered at $0$.
Moreover, $u$, which is also referred to
as the total field,  is the sum of a known incident field $u_i$
that is a solution of  the constant coefficient Helmholtz equation
\begin{equation}
\Delta u_i(x) + k^2 u_i(x) = 0 \  \ \mbox{for all}\  x \in \mathbb{R}^2
\label{introduction:incident}
\end{equation}
and an unknown scattered field $u_s$ which satisfies the Sommerfeld radiation condition
\begin{equation}
\lim_{r \to \infty}  \sup_{0 \leq t \leq 2\pi} \sqrt{r} \left|\frac{\partial u_s}{\partial r}
(r,t) - i k  u_s(r,t) \right| = 0.
\label{introduction:sommerfeld}
\end{equation}
In (\ref{introduction:sommerfeld}) and in what follows,
$r$ and $t$ are the usual polar coordinates for the point $x \in \mathbb{R}^2$
so that $x = r \exp(it)$.
Together (\ref{introduction:helm1}), (\ref{introduction:incident}) and (\ref{introduction:sommerfeld})
imply that the scattered field satisfies  the boundary value problem
\begin{equation}
\left\{
\begin{aligned}
&\Delta u_s(x) + k^2 \left(1+q(x)\right) u_s(x) = 
-k^2 q(x) u_i(x) \ \ \mbox{for all} \ \ x \in \mathbb{R}^2\\
&\sup_{0 \leq t \leq 2\pi}
\sqrt{r} \left|\frac{\partial u_s}{\partial r}
(r,t) - i k u_s(r,t) \right| = 0,
\end{aligned}
\right.
\label{introduction:scatbvp}
\end{equation}
and it is well known that (\ref{introduction:scatbvp}) suffices to uniquely
determine $u_s$ 
(see, for instance, Theorem~8.7 in Chapter~8 of \cite{Colton-Kress2}).

The total field $u$ and scattered field $u_s$ become
increasingly oscillatory as $k$ grows, and
 $\mathcal{O}\left(k\right)$ unknowns per spatial dimension are  required
to discretize them accurately using standard approaches, such as sampling their values
or expanding them in series of orthogonal polynomials.
Consequently, almost all numerical methods for 
the  solution of  (\ref{introduction:scatbvp}) have running times 
which grow at least quadratically with the wavenumber $k$.   
The principal exceptions are numerical-asymptotic methods, which 
use analytic information about the Helmholtz equation
to derive more efficient representations of its solutions which involve
fewer unknowns.  Some schemes of this type have running times which
are linear or even sublinear in the wavenumber $k$.
However, while many numerical-asymptotic schemes for solving the constant coefficient Helmholtz equation
have been proposed   (see, for instance, the survey
article \cite{Chandler-Wilde}), few apply in the case
of the variable coefficient Helmholtz equation.

Here, we describe a method for solving the variable coefficient Helmholtz
equation whose running time is $\mathcal{O}\left(k\log(k)\right)$ in typical cases 
but which only applies when the scattering potential  $q$ is radially symmetric.
Our approach exploits this symmetry to construct an expansion of the total field  through  separation of 
variables.  More explicitly,  we represent $u$ in the interior
of the the disk $\Omega$ via a sum of the form
%
\begin{equation}
u(r,t) = \sum_{n=-m}^m a_n \psi_{\left|n\right|}(r) \exp(int),
\label{introduction:totalrep}
\end{equation}
where, for each nonnegative integer $n$, $\psi_n$ satisfies a second order differential equation 
whose coefficients depend on the wavenumber $k$ and the index $n$ and which we refer
to as the ``perturbed Bessel equation.''
The value of $m$ in (\ref{introduction:totalrep}) is equal to the number of Fourier modes
needed to accurately represent the restriction of the incident wave $u_i$ to the boundary
of the scatterer $\partial\Omega$.  Since $u_i$ satisfies the constant coefficient
Helmholtz equation at wavenumber $k$,  $m$ is  $\mathcal{O}\left(k\right)$ in typical applications.
Moreover, we solve the perturbed Bessel equation
using a method whose worst case running time  appears to be  $\mathcal{O}\left(\log(k)\right)$.

We say  ``appears to be''  because our evidence for this claim is experimental.
To construct the $\psi_n$, we use an approach based on the method of \cite{BremerKummer} 
for the numerical solution of second order differential equations of
the form
\begin{equation}
y''(t) +  k^2 \eta(t) y(t) = 0 \ \ \mbox{for all} \ \ a < t < b.
\label{introduction:ode}
\end{equation}
When $\eta$ is positive,  the solutions of (\ref{introduction:ode}) are oscillatory and 
$\mathcal{O}\left(k\right)$
unknowns are required to discretize them using standard approaches (e.g., by sampling their
values).  The solver of \cite{BremerKummer} instead
represents  them using
 a nonoscillatory phase function which can be calculated and stored efficiently.
Indeed, according to the estimates of \cite{Bremer-Rokhlin}, 
the method of \cite{BremerKummer}  runs in time independent of $k$
\emph{under the assumption that $\eta$ is smooth and strictly positive on the interval $[a,b]$}.
The equation which defines $\psi_n$ can have turning points, however, and in the  neighborhood
of such a point, the estimates of \cite{Bremer-Rokhlin} do not apply.
Here, we present the results of numerical experiments which strongly indicate that at worst
our method requires  $\mathcal{O}\left(\log(k)\right)$ to construct    each $\psi_n$.
Assuming this is correct, $\mathcal{O}\left(k \log(k)\right)$ time is required to construct
all of the functions $\psi_n$ appearing in (\ref{introduction:totalrep}).
A further $\mathcal{O}\left(k \log(k)\right)$ operations are needed to calculate the coefficients
in this expansion, so that the time required to solve
(\ref{introduction:scatbvp}) is $\mathcal{O}\left(k \log(k)\right)$.

It is well known that the condition number of evaluation of oscillatory functions increases
with the frequency of their oscillations and that this generally limits the relative accuracy 
with which they can evaluated numerically (see, for instance, \cite{Higham}). 
As a consequence, numerical schemes for solving (\ref{introduction:scatbvp}),
including the algorithm of this paper, loses accuracy with increasing $k$.
Section~\ref{section:experiments} of this article describes
numerical experiments which were conducted to assess the speed and accuracy of our scheme.



Separation of variables is  hardly a new idea; but, it appears to have been
rarely used as a numerical tool  for solving (\ref{introduction:scatbvp}).
This is most likely because it offers no advantage in asymptotic 
running time when standard numerical methods are used to 
solve  the differential equations which arise.  Moreover, when
the parameter $n$ is  large, 
 $\psi_n$ behaves like a rapidly increasing
exponential function on some part of its domain, and the accurate numerical
solution of the corresponding differential equation becomes 
difficult. One of the few robust numerical algorithms for solving (\ref{introduction:scatbvp})
via separation of variables is described in   \cite{Hoskins-Rokhlin}.
It proceeds by reformulating the differential equations that arise 
an integral equations which are then inverted using a ``fast direct solver.''
This results in an extremely accurate and robust method
for solving (\ref{introduction:scatbvp}); however,
the running time of this method still grows quadratically with $k$.
Our scheme also bears some similarities to that of \cite{Gillman-Barnett-Martinsson}.   There,
spectral methods are used to construct expansions of the total field in the interior
of $\Omega$ and of the scattered field in the exterior of $\Omega$, and the coefficients
in these expansion are found using a mechanism very similar to that used by our
algorithm to compute the coefficients in the expansion  (\ref{introduction:totalrep}).  However,
the scheme of \cite{Gillman-Barnett-Martinsson} is not limited to radially symmetric potentials
and its running time grows somewhat faster than $\mathcal{O}\left(k^2\right)$.

The remainder of this document is structured as follows.
Section~\ref{section:determination} discusses the solution of 
of the boundary value problem (\ref{introduction:scatbvp}) via the method separation of variables.
Our algorithm for the numerical solution of 
the boundary value problem (\ref{introduction:scatbvp}) is discussed in
Sections~\ref{section:algorithm} and \ref{section:pbessel}.
In Section~\ref{section:experiments}. we describe several numerical experiments which
demonstrate the properties of our solver. We close in Section~\ref{section:conclusion}
with a few brief remarks regarding this work and 
a discussion of possible directions for future research.

\label{section:introduction}
\end{section}

\begin{section}{Determination of the Scattered Field Through Separation of Variables}
\label{section:determination}


The total field satisfies
\begin{equation}
\Delta u(r,t)  + k^2 ( 1 + q(r) ) u(r,t) = 0
\label{determination:toteq}
\end{equation}
in $\Omega$ and is nonsingular at $0$.  Separating variables 
in (\ref{determination:toteq}) gives us the representation formula
\begin{equation}
u(r,t) = \sum_{n-\infty}^\infty a_n \psi_{\left|n\right|}(r) \exp(int),
\label{determination:totrep}
\end{equation}
where, for each nonnegative integer $n$, $\psi_n$ is a nonzero solution of
\begin{equation}
r^2 y''(r) + r y'(r) + \left( k^2( 1 + q(r)) r^2 - n^2 \right) y(r) = 0,
\ \ 0 < r < R,
\label{determination:psieq}
\end{equation}
which is regular at $0$.
When $q=0$, (\ref{determination:psieq}) becomes Bessel's differential equation
and we refer to it as the perturbed Bessel equation.
Since this equation has a regular singular point at 0 and the roots of the indicial equation
are $n$ and $-n$, there exists a 
basis $\{u_n,v_n\}$  in the space of its solutions such that
\begin{equation}
u_n(r) = \mathcal{O}\left(r^n\right) \ \ \mbox{as}\ \ r \to 0 
\label{determination:psi0}
\end{equation}
and
\begin{equation}
v_n(r) = \begin{cases}
\mathcal{O}\left(\log(r)\right) \ \ \mbox{as} \ \ r \to 0 & \mbox{if} \ n = 0 \\
\mathcal{O}\left(r^{-n}\right)   \ \ \ \ \ \mbox{as} \ \ r \to 0  & \mbox{if} \ n > 0
\end{cases}
\end{equation}
(see, for instance, Chapter~5 of \cite{Hille}).
Since  $\psi_n$, is nonsingular at $0$, 
it must be a multiple of the solution $u_n$.   Because we do not impose
a second boundary condition on $\psi_n$, it is only determined
up to a (nonzero) multiplicative constant.  
It is perhaps tempting to impose a second condition on $\psi_n$ in the hopes
of uniquely determining it.  However, doing so generally leads to a boundary value
problem which is not solvable for all possible $k$.
Moreover, this ambiguity has no
impact on the form of the expansion (\ref{determination:totrep}) used to
represent the total field in the interior of $\Omega$, although
the particular value of the coefficient $a_n$ depends on the
choice of $\psi_n$.

Because $q$ is supported inside of the disk $\Omega$, the
scattered field $u_s$ solves
\begin{equation}
\left\{
\begin{aligned}
&\Delta u_s(r,t) + k^2 u_s(r,t) =  0  
\ \ \mbox{in} \ \ \Omega^c\\
&\lim_{r \to \infty} 
\sup_{0 \leq t \leq 2\pi}
\sqrt{r} \left|\frac{\partial u_s}{\partial r}
(r,t) - i  k u_s(r,t) \right| = 0.
\end{aligned}
\right.
\label{determination:scatbvp}
\end{equation}
Separating variables in (\ref{determination:scatbvp})
shows that $u_s$ can be represented in the exterior of $\Omega$
via a sum of the form
\begin{equation}
u_s(r,t) = \sum_{n=-\infty}^\infty b_n \sigma_{|n|}(k r)\exp(int),
\label{determination:scatrep0}
\end{equation}
where, for each nonnegative $n$, $\sigma_n$ is a solution of Bessel's differential
\begin{equation}
r^2 y''(r) + r y'(r) + \left(k^2r^2 - n^2\right)y(r) = 0, \ \ R < r < \infty,
\label{determination:bessel}
\end{equation}
consistent with the Sommerfeld radiation condition.  It follows from the 
integral representation formula
\begin{equation}
H_n(z) = \exp(iz) \left( \frac{-2^{n+1} i}{\sqrt{\pi}\ \Gamma\left(n+\frac{1}{2}\right)}     
z^n \int_0^\infty (x^2-ix)^{n-\frac{1}{2}} \exp(-2xz)\ dx\right)
\label{determination:hn}
\end{equation}
for the Hankel function of the first kind of order $n$ 
(which  can be found as Formula~3.388(4) in \cite{Gradshteyn})
and standard results regarding 
the decay properties of Laplace transforms (see, for instance, \cite{Widder}) that 
$H_n(kr)$  is a solution of (\ref{determination:bessel})
which has the appropriate behavior at infinity.
In particular, $u_s$ admits the representation
\begin{equation}
u_s(r,t) = \sum_{n=-\infty}^\infty b_n H_n(k r)\exp(int)
\label{determination:scatrep}
\end{equation}
in the exterior of $\Omega$.

Standard elliptic regularity results (see, for instance, Section~8.4 of \cite{Gilbarg-Trudinger})
imply that   the total field $u$ and its derivative
$\partial u / \partial r$ with respect to the radial variable
are continuous   across the boundary  $\partial \Omega$ of the disk $\Omega$.
This suffices to determine the coefficients
 $a_n$ in (\ref{determination:totrep}) and $b_n$  in (\ref{determination:scatrep}).
To see this, we first let
\begin{equation}
u_i(R,t) = \sum_{n=-\infty}^\infty c_n \exp(int)
\label{determination:ui}
\end{equation}
and
\begin{equation}
\frac{\partial u_i}{ \partial r } (R,t) \approx \sum_{n=-\infty}^\infty d_n \exp(int)
\label{determination:uiprime}
\end{equation}
be the Fourier expansions of the restrictions of the incident field $u_i$
and its radial derivative to $\partial\Omega$.
Since $u - u_s  = u_i$,  the continuity of $u$ and its radial derivative
imply that
\begin{equation}
\sum_{n=-\infty}^\infty a_n \psi_{|n|}(R) \exp(int) - 
\sum_{n=-\infty}^\infty b_n H_{n}(k R) \exp(int) =
\sum_{n=-\infty}^\infty c_n \exp(int)
\label{determination:eq1}
\end{equation}
and
\begin{equation}
\sum_{n=-\infty}^\infty a_n \psi_{|n|}'(R) \exp(int) - 
\sum_{n=-\infty}^\infty b_n k H_{n}'(k R) \exp(int) =
\sum_{n=-\infty}^\infty d_n \exp(int)
\label{determination:eq2}
\end{equation}
for all $0 \leq t \leq 2\pi$.    Owing to the orthogonality of the set  $\{\exp(int)\}$,
(\ref{determination:eq1}) and (\ref{determination:eq2}) hold if and only if
for each integer $n$, $a_n$ and  $b_n$ satisfy the linear system of equations
\begin{equation}
\left\{
\begin{aligned}
a_n \psi_{|n|}(R) - b_n H_n\left(k R\right) &= c_n \\
a_n \psi_{|n|}'(R) - b_n k H_n'\left(k R\right) &= d_n. \\
\end{aligned}
\right.
\label{determination:sys}
\end{equation}
Since the scattered field is uniquely determined by (\ref{introduction:scatbvp}),
each of the systems  (\ref{determination:sys}) must be uniquely solvable.
In particular, for each integer $n$, the determinant 
\begin{equation}
H_n\left(k R\right) \psi_{|n|}'(R) - k\ \psi_{|n|}(R) H_n'\left(k R\right)
\end{equation}
of the coefficient matrix in (\ref{determination:sys}) is necessarily nonzero
and the coefficients $a_n$ and $b_n$ are given by the formulas
\begin{equation}
a_n = 
\frac
{-k H_n'\left(k R\right) c_n +  H_n\left(k R\right) d_n}
{H_n\left(k R\right) \psi_{|n|}'(R) - k\ \psi_{|n|}(R) H_n'\left(k R\right) }\\
\label{determination:an}
\end{equation}
and
\begin{equation}
b_n = 
\frac
{-\psi_{|n|}'(R) c_n +  \psi_{|n|}(R) d_n}
{H_n\left(k R\right) \psi_{|n|}'(R) - k\ \psi_{|n|}(R) H_n'\left(k R\right) }.
\label{determination:bn}
\end{equation}

This procedure determines the scattered field in the sense that
once the coefficients $a_n$ and $b_n$ are determined through (\ref{determination:an})
and (\ref{determination:bn}), the scattered field can be evaluated
at any point outside of $\Omega$ via (\ref{determination:scatrep})
and it can be evaluated at any point inside of $\Omega$ or on $\partial\Omega$
by first evaluating the total field via  (\ref{determination:totrep})
and then subtracting the value of the (known) incident field $u_i$
from the result.

\end{section}

\begin{section}{A Numerical Algorithm for the Determination of the Scattered Field}
\label{section:algorithm}

Our algorithm operates in two phases: a precomputation phase in which the 
perturbed Bessel equation is repeatedly solved in order
to construct the functions $\psi_n$ used in the representation
of the total field, and a solution phase
in which the coefficients in the expansions of the total and scattered fields
are calculated.  We refer to the first procedure as the ``precomputation phase'' because
in many applications it is necessary to solve (\ref{introduction:scatbvp}) for multiple
incident fields $u_i$ while the wavenumber $k$ and scattering potential $q$
are fixed.  In this event, the precomputation phase is only executed once 
and the solution phase is  executed once for each incident field.

The precomputation phase takes as input the wavenumber $k$ for the 
problem, a subroutine for evaluating the scattering potential
$q$ at any specified point, a list 
\begin{equation}
\chi_1 < \chi_2 < \ldots < \chi_s
\label{algorithm:singpts}
\end{equation}
of all  of the  points  on the interval $(0,R)$ at which the scattering potential $q(r)$ 
is nonsmooth (recall that $q$ is assumed to be piecewise smooth),
and a positive integer $m$.    It consists of calculating functions
\begin{equation}
\psi_0,\ \psi_1,\ \ldots,\ \psi_m
\end{equation}
such that for each  integer $0 \leq n  \leq m$, $\psi_n$ is a solution
of the perturbed Bessel equation (\ref{determination:psieq}) 
which is regular at $0$.  We detail our algorithm for solving
the perturbed Bessel equation in Section~\ref{section:pbessel}, which
follows this one.  Based on strong experimental evidence, we believe that the running time
of the precomputation phase is $\mathcal{O}\left(m\log(m)\right)$.

The solution phase takes as input a routine for evaluating
the  incident field $u_i$ and its derivative 
$\partial u_i / \partial r$ with respect to the radial variable $r$.   
It proceeds by first forming the approximations
\begin{equation}
u_i(R,t) \approx \sum_{n=-m}^m c_n \exp(int)
\label{algorithm:ui}
\end{equation}
and
\begin{equation}
\frac{\partial u_i}{ \partial r } (R,t) \approx \sum_{n=-m}^m  d_n \exp(int)
\label{algorithm:uider}
\end{equation}
of the restrictions of $u_i$ and its normal derivative to the boundary
$\partial\Omega$ in the usual way ---   that is, using the fast Fourier transform.
The integer $m$ must be sufficiently large for the 
approximations (\ref{algorithm:ui}) and (\ref{algorithm:uider}) to be highly accurate.
In the next step of the solution phase, the coefficients in the truncated expansion
\begin{equation}
u(r,t) = \sum_{n=-m}^m a_n \psi_{\left|n\right|}(r) \exp(int)
\label{algorithm:utot}
\end{equation}
which represents the total wave in the interior of $\Omega$
are computed using Formula~(\ref{determination:an}).
Finally, the coefficients in the truncated expansion
\begin{equation}
u_s(r,t) = \sum_{n=-m}^m b_n H_n(k r)\exp(int)
\label{algorithm:uscat}
\end{equation}
used to represent  the scattered field in the exterior of  $\Omega$ 
are computed using (\ref{determination:bn}).  
The fast Fourier transforms  take $\mathcal{O}\left(m \log(m)\right)$ operations, and 
they  dominate the cost of this phase of the algorithm.
  Only $\mathcal{O}\left(m\right)$
operations are required to construct the coefficients in the expansions (\ref{algorithm:utot}) and 
(\ref{algorithm:uscat}).

The coefficients in the expansions (\ref{algorithm:utot}) and  (\ref{algorithm:uscat})
are the principal outputs of the algorithm of this paper.  Once they have been determined,
 the scattered field can be evaluated at any point in the exterior of $\Omega$
in $\mathcal{O}\left(m\right)$ operation  by evaluating the sum (\ref{algorithm:uscat}),
and it can be evaluated at any point in the interior of $\Omega$ in $\mathcal{O}\left(m\right)$
operations by evaluating (\ref{algorithm:utot}) and subtracting the value of the incident field
$u_i$.

Since $u_i$ is a solution of the constant coefficient Helmholtz
equation at wavenumber $k$, it is expected that $m$ will be on the order of $k$
so that the running time of our algorithm is $\mathcal{O}\left(k\log(k)\right)$.
We found $m= \pi/2 R k$ to be sufficient in all of the numerical experiments discussed in this paper.
If $m$ is not known {\it a priori} it can be determined through an adaptive
procedure.   For instance, starting from an initial guess, $m$ could be gradually
increased until the coefficients in the expansions
(\ref{algorithm:ui}) and  (\ref{algorithm:uider}) decay sufficiently fast.

\end{section}

\begin{section}{Numerical Solution of the Perturbed Bessel Equation}
\label{section:pbessel}

It can be easily seen  that if $\psi_n$ solves (\ref{determination:psieq}),
then $\varphi_n(r) = \sqrt{r}\ \psi_n(r)$ is a solution of
\begin{equation}
y''(r) + Q(r) y(r) = 0, \ \ 0 < r < R,
\label{pbessel:ode}
\end{equation}
where
\begin{equation}
Q(r) = \lambda^2 \left( 1 + q(r) \right) + \frac{\frac{1}{4}-n^2}{r^2}.
\label{pbessel:coef}
\end{equation}
We refer to (\ref{pbessel:ode}) as the normal form of the perturbed Bessel equation
and we find it more convenient to work with than (\ref{determination:psieq}).
Among other things, a great deal of information about the behavior of the solutions
of (\ref{pbessel:ode}) can be easily discerned from the coefficient (\ref{pbessel:coef}).
Indeed, according to standard asymptotic results (see, for instance, \cite{Olver} or \cite{Fedoryuk}),
in intervals on which $Q$ is positive the solutions of (\ref{pbessel:ode})
behave roughly as oscillatory exponential functions, while they resemble increasing or 
decreasing exponential functions in  intervals on which $Q$ is negative.
The zeros of $Q$ which separate these regions are known as turning points for (\ref{pbessel:ode}).
We do not, in fact, produce a solution of the  perturbed Bessel equation
over the entire interval $[0,R]$ since $Q$ is singular at $0$ 
and, as a consequence, all but one of the solutions of (\ref{pbessel:ode}) are also singular at $0$.
Instead, we produce a solution over $[10^{-15},R]$,  which generally 
suffices for the purposes of  numerical computation.

The first step of our algorithm for solving
the perturbed Bessel equation consists of forming a  partition
\begin{equation}
10^{-15}= \xi_1 < \xi_2 < \ldots < \xi_t = R
\end{equation}
of the solution interval $[10^{-15},R]$ such that $Q$ is smooth and does not change sign on each interval
$(\xi_j,\xi_{j+1})$.  We do so by finding the set of all zeros of $Q$ on $[0,R]$ 
and merging it with the list (\ref{algorithm:singpts}) of the singularities
of $Q$ provided by the user as well as the additional points $10^{-15}$ and $R$.

Next, for each $j=1,\ldots,t$, our solver forms a basis 
\begin{equation}
\mathscr{B}_j = \left\{u_j, v_j\right\}
\end{equation}
in space of the restrictions of solutions of (\ref{pbessel:ode}) to $\left(\xi_j,\xi_{j+1}\right)$.
The mechanisms used to construct and represent these basis functions differ depending on whether
$Q$ is positive in the interval (the oscillatory regime) or negative in the 
interval (the nonoscillatory regime).  We discuss the details in each case below.

Finally, for each $j=1,\ldots,t$, it calculates coefficients $\gamma_j$ and $\zeta_j$ such that
the restriction of the desired solution $\varphi_n$ of (\ref{pbessel:ode})
to the interval $\left(\xi_j,\xi_{j+1}\right)$ is 
\begin{equation}
\gamma_j u_j(r) + \eta_j v_j(r).
\end{equation}
In the case of the first interval $\left(\xi_1, \xi_2\right)$,
the values of $\varphi_n$ and its derivative at the point $\xi_1 = 10^{-15}$ are first estimated
via  the asymptotic approximation
\begin{equation}
\varphi_n(r) \sim \sqrt{r} J_n \left(\sqrt{q(0) + k^2} r \right) \ \ \mbox{as} \ \ r \to 0,
\end{equation}
which can be easily derived using standard methods (see, for instance, \cite{Olver} or \cite{Fedoryuk}).
Then,  the linear system of equations
\begin{equation}
\left\{
\begin{aligned}
\gamma_{1} u_{1}(\xi_1) + \zeta_{j-1} u_{1} ( \xi_1) =  \varphi_n(\xi_1)  \\
\gamma_{1} u_{1}'(\xi_1) + \zeta_{j-1} u_{1}' ( \xi_1) =  \varphi_n'(\xi_1)
\end{aligned}
\right.
\end{equation}
is solved for $\gamma_1$ and $\zeta_1$.
For each $j >1$, the coefficients $\gamma_j$ and $\zeta_j$
are determined by enforcing the continuity of $\varphi_n$ and its derivative at the point $\xi_{j}$.  
More explicitly,  we solve the system of linear equations
\begin{equation}
\left\{
\begin{aligned}
\gamma_{j-1} u_{j-1}(\xi_j) + \eta_{j-1} u_{j-1} ( \xi_j) =  \gamma_{j} u_{j}(\xi_j) + \eta_{j} u_{j} ( \xi_j)  \\
\gamma_{j-1} u_{j-1}'(\xi_j) + \eta_{j-1} u_{j-1}' ( \xi_j) =  \gamma_{j} u_{j}'(\xi_j) + \eta_{j} u_{j}' ( \xi_j)  \\
\end{aligned}
\right.
\end{equation}
for $\gamma_j$ and $\eta_j$.

\begin{subsection}{The oscillatory regime}

For intervals $\left(\xi_j,\xi_{j+1}\right)$ on which 
$Q$ is positive, we use a basis generated by a nonoscillatory phase function.
A function $\alpha$ is a phase function for  the differential equation
\begin{equation}
y''(t) + \eta(t) y(t) = 0 \ \ \mbox{for all} \ \ a < t <b
\label{pbessel:2nd}
\end{equation}
provided $\alpha'(r) > 0$ for all $a < r <b$ and 
\begin{equation}
\frac{\sin\left(\alpha(r)\right)}{\sqrt{\alpha'(r)}} \ \ \mbox{and} \ \
\frac{\cos\left(\alpha(r)\right)}{\sqrt{\alpha'(r)}}
\end{equation}
is a basis in its space of solutions.       An extensive discussion of 
phase functions for second order differential equations
can be found in \cite{Neuman}.     Among other things, it is shown there that $\alpha$ is a phase function
for (\ref{pbessel:2nd}) if and only if its derivative satisfies the 
nonlinear second order differential equation
\begin{equation}
\left(\alpha'(r)\right)^2 = \eta(r) - \frac{1}{2}\frac{\alpha'''(r)}{\alpha'(r)}
+ \frac{3}{4} \left(\frac{\alpha''(r)}{\alpha'(r)}\right)^2,\ \ \ a <  r < b,
\label{pbessel:kummer}
\end{equation}
which we refer to as Kummer's equation after E.~E.~Kummer who studied
it in \cite{Kummer}.     We note that (\ref{pbessel:kummer}) only determines
$\alpha$ up to a constant.  For our purposes, the constant is largely irrelevant
(our only requirement is that it not be too large in magnitude),
and  we always determine it by taking $\alpha$ to be zero at left-hand
endpoint of the interval on which it is defined.

In \cite{Bremer-Rokhlin}, it is shown
that, under mild assumptions on the coefficient $\eta$ (including the condition
that it be positive on the interval $[a,b]$), there exists
 a phase function for (\ref{pbessel:2nd}) which is roughly as oscillatory
as the coefficient $\eta$.  
Moreover, in \cite{BremerKummer}, a fast algorithm
for the numerical calculation of this nonoscillatory phase functions
is presented.  It operates by first introducing a ``windowed version''
$\tilde{\eta}$ of $\eta$ such that
\begin{equation}
\tilde{\eta}(r)  = 
\begin{cases}
\lambda^2 & \mbox{for all} \ \ a < r < \frac{3a+b}{4}  \\
\eta(r) & \mbox{for all} \ \ \frac{a+3b}{4} < r < b
\end{cases}
\end{equation}
with $\lambda$ a constant chosen to be roughly on the order of $\sqrt{\eta(a)}$.
Since $\tilde{\eta}=\lambda^2$ near $a$, the nonoscillatory
phase function for the equation
\begin{equation}
y''(r) + \tilde{\eta}(r)y(r) = 0
\end{equation}
is equal to  $\lambda r$ near $a$.
By solving the initial value problem
%
\begin{equation}
\left\{
\begin{aligned}
\left(\tilde{\alpha}'(r)\right)^2 &= \tilde{\eta}(r) - \frac{1}{2}\frac{\tilde{\alpha}'''(r)}{\tilde{\alpha}'(r)}
+ \frac{3}{4} \left(\frac{\tilde{\alpha}''(r)}{\tilde{\alpha}'(r)}\right)^2\\
\tilde{\alpha}'(r) &= \lambda \ \ \mbox{and}\ \ \tilde{\alpha}''(r) = 0,
\end{aligned}
\right.
\label{pbessel:initial}
\end{equation}
the values of $\tilde{\alpha}'(b)$ and $\tilde{\alpha}''(b)$
are determined.  Since $\tilde{\eta}$ is equal to $\eta$ near $b$, these
values closely approximate $\alpha'(b)$ and $\alpha''(b)$.
The function $\alpha$ is then determined over the interval $(a,b)$ by solving
the terminal value problem
\begin{equation}
\left\{
\begin{aligned}
\left(\alpha'(r)\right)^2 &= \eta(r) - \frac{1}{2}\frac{\alpha'''(r)}{\alpha'(r)}
+ \frac{3}{4} \left(\frac{\alpha''(r)}{\alpha'(r)}\right)^2\\
\alpha'(b) &= \tilde{\alpha}'(b)\\
 \alpha''(b) &= \tilde{\alpha}''(b).
\end{aligned}
\\
\right. 
\label{pbessel:terminal}
\end{equation}
To solve (\ref{pbessel:initial}) and (\ref{pbessel:terminal}), we 
use a spectral method which represents the phase function $\alpha$
and its derivatives using piecewise Chebyshev expansions over a collection of
subintervals of $(a,b)$.  The subintervals are chosen adaptively.  Other mechanisms
for the solution of these differential equations could be used, so long as 
they are well-suited for stiff problems.
When the coefficient $\eta$ is strictly positive on $[a,b]$, the nonoscillatory
phase function for (\ref{pbessel:2nd}) produced by this algorithm can
be constructed and evaluated in time independent of the magnitude of $\eta$ (which is a measure
of the frequency of oscillation of the solutions of (\ref{pbessel:2nd})).

Since (\ref{pbessel:ode}) can have turning points at the endpoints $\xi_j$ and
$\xi_{j+1}$, the estimates of \cite{Bremer-Rokhlin} do not apply.  However,
we have found experimentally  (see Sections~\ref{section:experiments:pbessel1}
and \ref{section:experiments:pbessel2})
that the worst case running time of the algorithm of \cite{BremerKummer}
is $\mathcal{O}\left(\log(k)\right)$ in this case.

\end{subsection}


\begin{subsection}{The nonoscillatory regime}

For intervals $(\xi_j,\xi_{j+1})$ on which $Q$ is negative, 
we use a basis  $\{u_j,v_j\}$  of solutions of  (\ref{pbessel:ode}) 
such that $u_j$ resembles an increasing exponential function and $v_j$ resembles
a decreasing exponential function.
Because the cost of representing these functions using standard methods, such as
through expansions
in orthogonal polynomials or via sampling their values, increases rapidly with $n$ and $k$,
we instead construct their logarithms and use these to evaluate
$u_j$ and $v_j$ as needed.

If $y(r) = \exp(\sigma(t))$ satisfies the second order differential equation
(\ref{pbessel:ode}), then it can be easily verified that  
$\sigma'$ satisfies the Riccati equation
\begin{equation}
\sigma''(r) + (\sigma'(r))^2 + Q(r) = 0
\label{pbessel:riccati}
\end{equation}
(see, for instance, Chapter~4 of \cite{Hille} for a discussion
of the Riccati equation).    To construct the logarithm of $u_j$, we first calculate a solution
$\sigma'$  of (\ref{pbessel:riccati})  which satisfies the initial condition
 $\sigma'\left(\xi_j\right) =0$.   Again, we use an adaptive
spectral solver whose output is a piecewise Chebyshev expansion representing the solution.
We then use spectral integration to form the antiderivative $\sigma$ of $\sigma'$
such that $\sigma(\xi_{j+1}) = 0$.  Since the dominant solution of 
(\ref{pbessel:riccati}) when solving in the forward direction is increasing,
 the function $u_j$ constructed in this fashion resembles an increasing exponential function.

To construct the logarithm of $v_j$, we solve (\ref{pbessel:riccati}), imposing the terminal
condition $\sigma'\left(\xi_{j+1}\right) = 0$.
Next we use spectral integration to form the antiderivative $\sigma$ of $\sigma'$
such that $\sigma(\xi_{j}) = 0$.  Since the dominant solution of (\ref{pbessel:riccati})
in the backward direction is decreasing, $v_j$ resembles a decreasing exponential
function.

Based on the extensive numerical experiments of Section~\ref{section:experiments},
we believe that the worst case running time of this procedure for 
constructing $u_j$ and $v_j$ is $\mathcal{O}\left(\log(k)\right)$.


\end{subsection}

\end{section}

\begin{figure}[p]
\small

\hfill
\includegraphics[width=.40\textwidth]{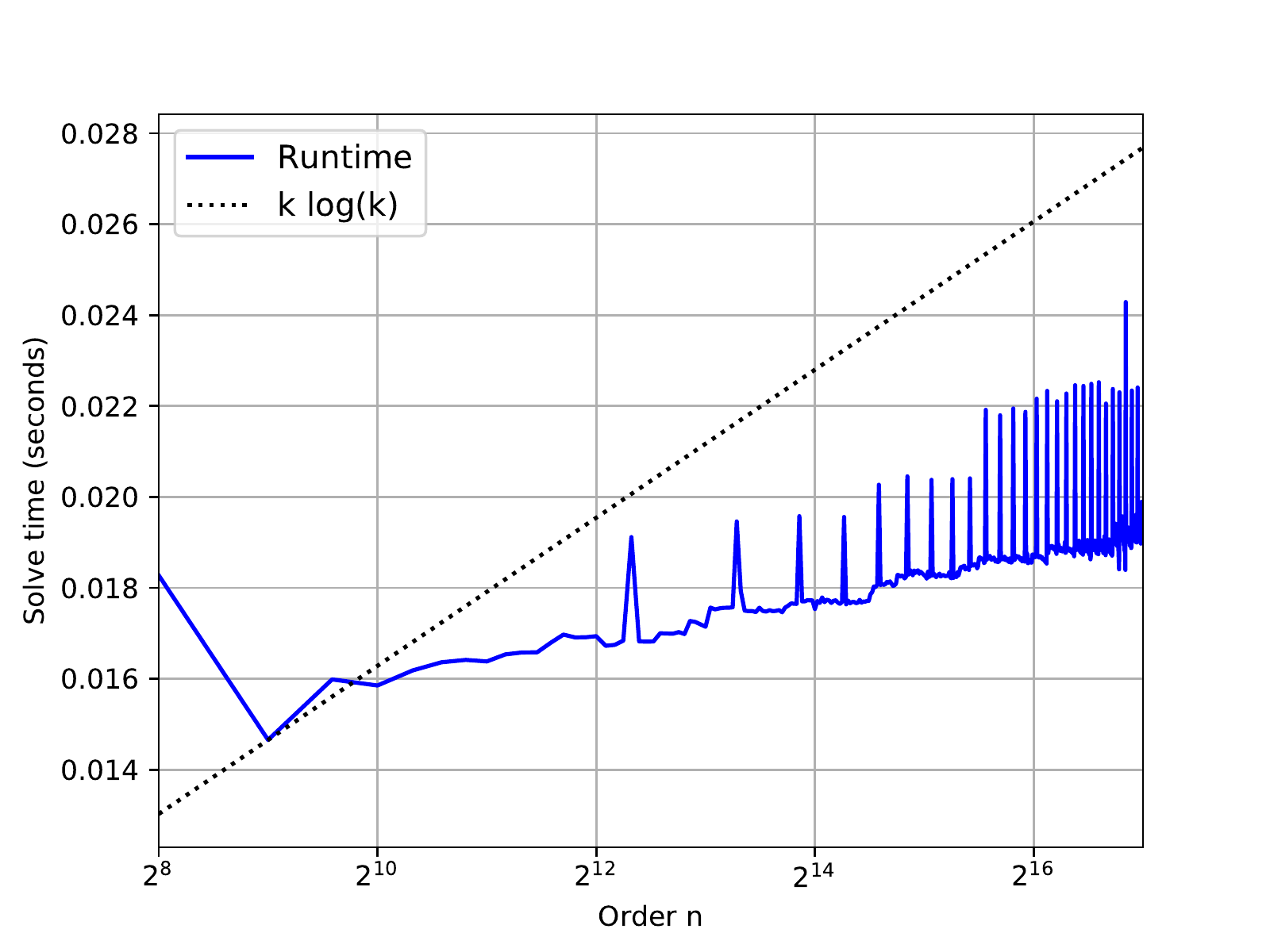}
\hfill
\includegraphics[width=.40\textwidth]{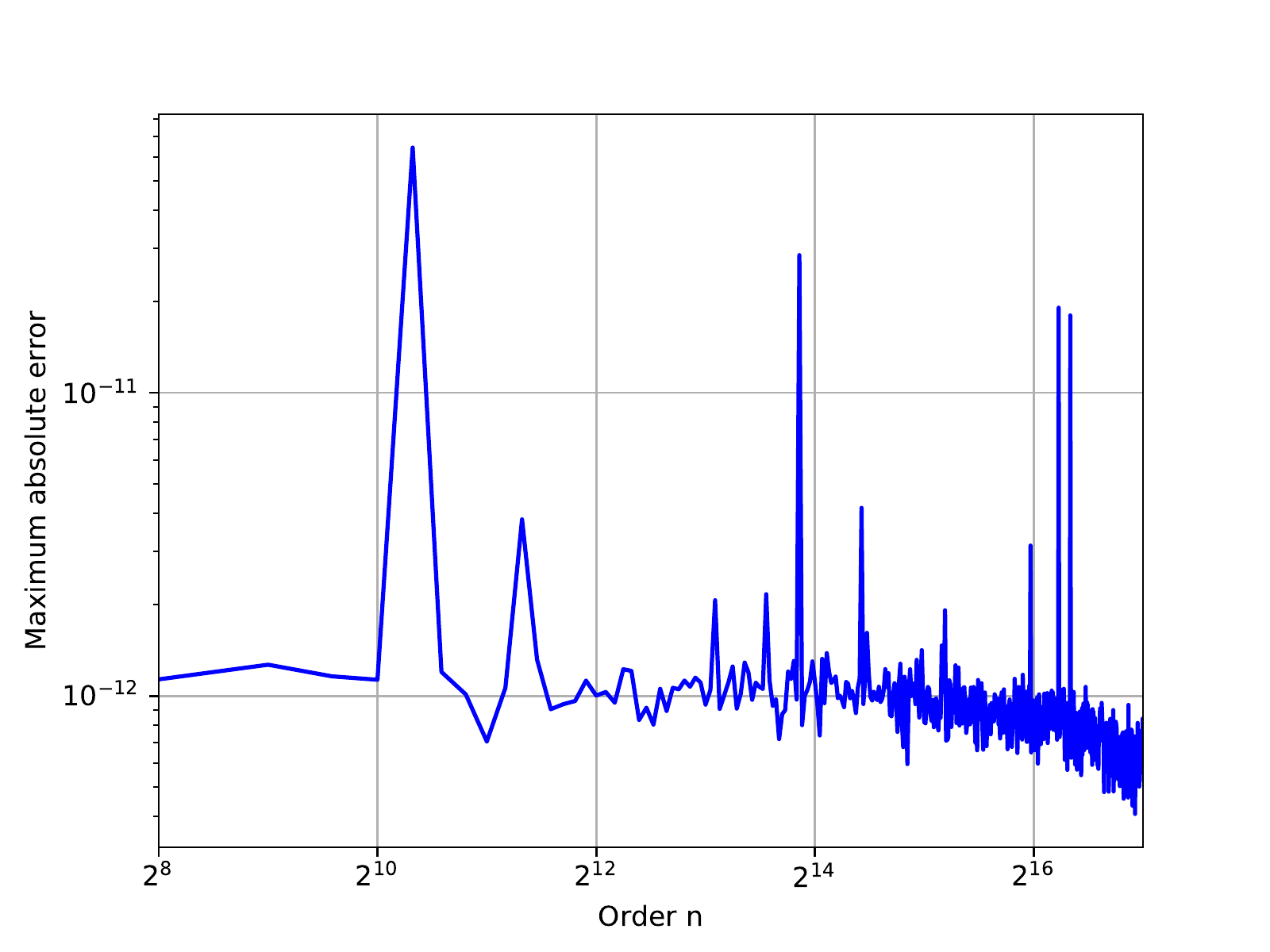}
\hfill

\hfill
\includegraphics[width=.40\textwidth]{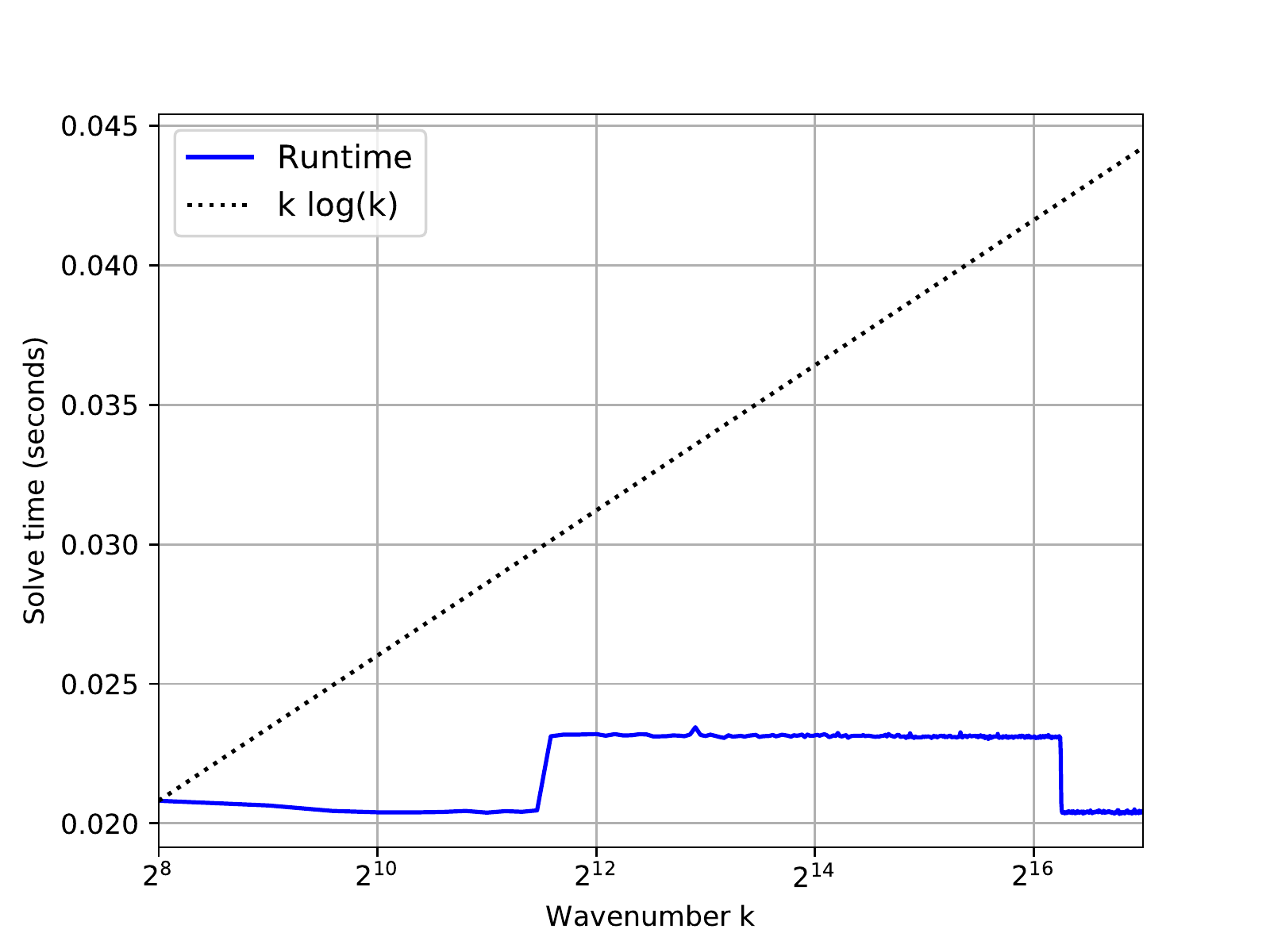}
\hfill
\includegraphics[width=.40\textwidth]{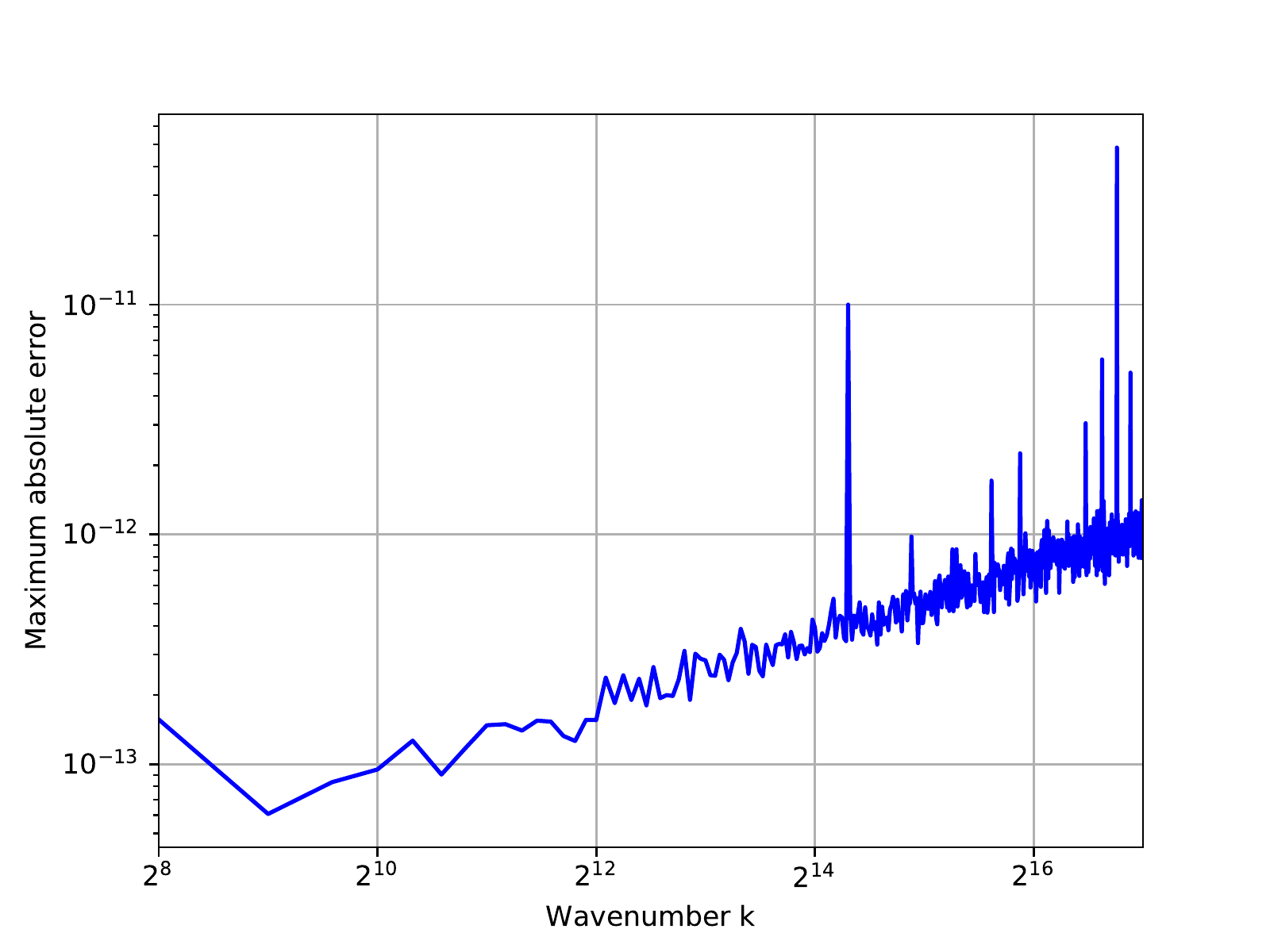}
\hfill

\hfill
\includegraphics[width=.40\textwidth]{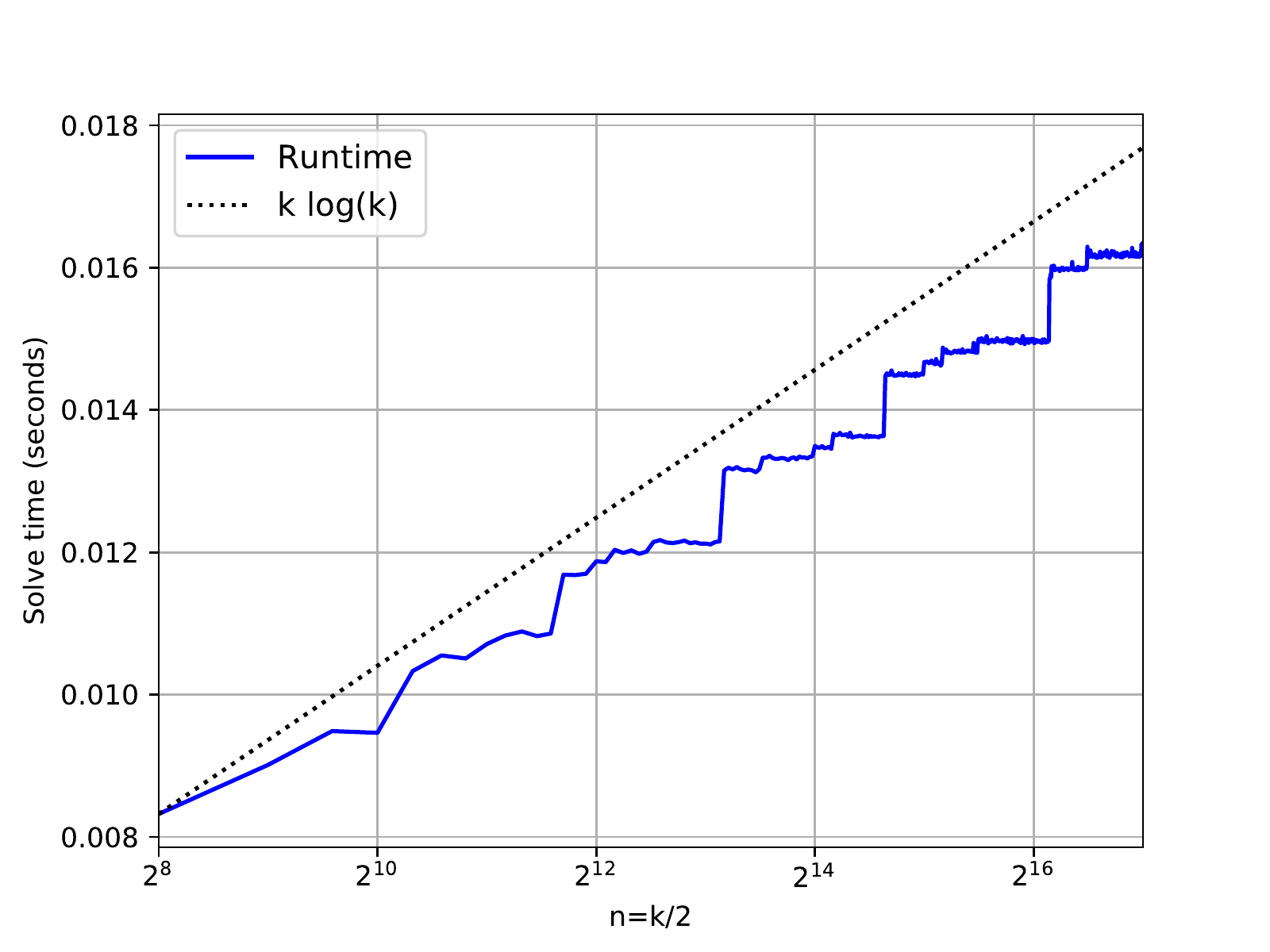}
\hfill
\includegraphics[width=.40\textwidth]{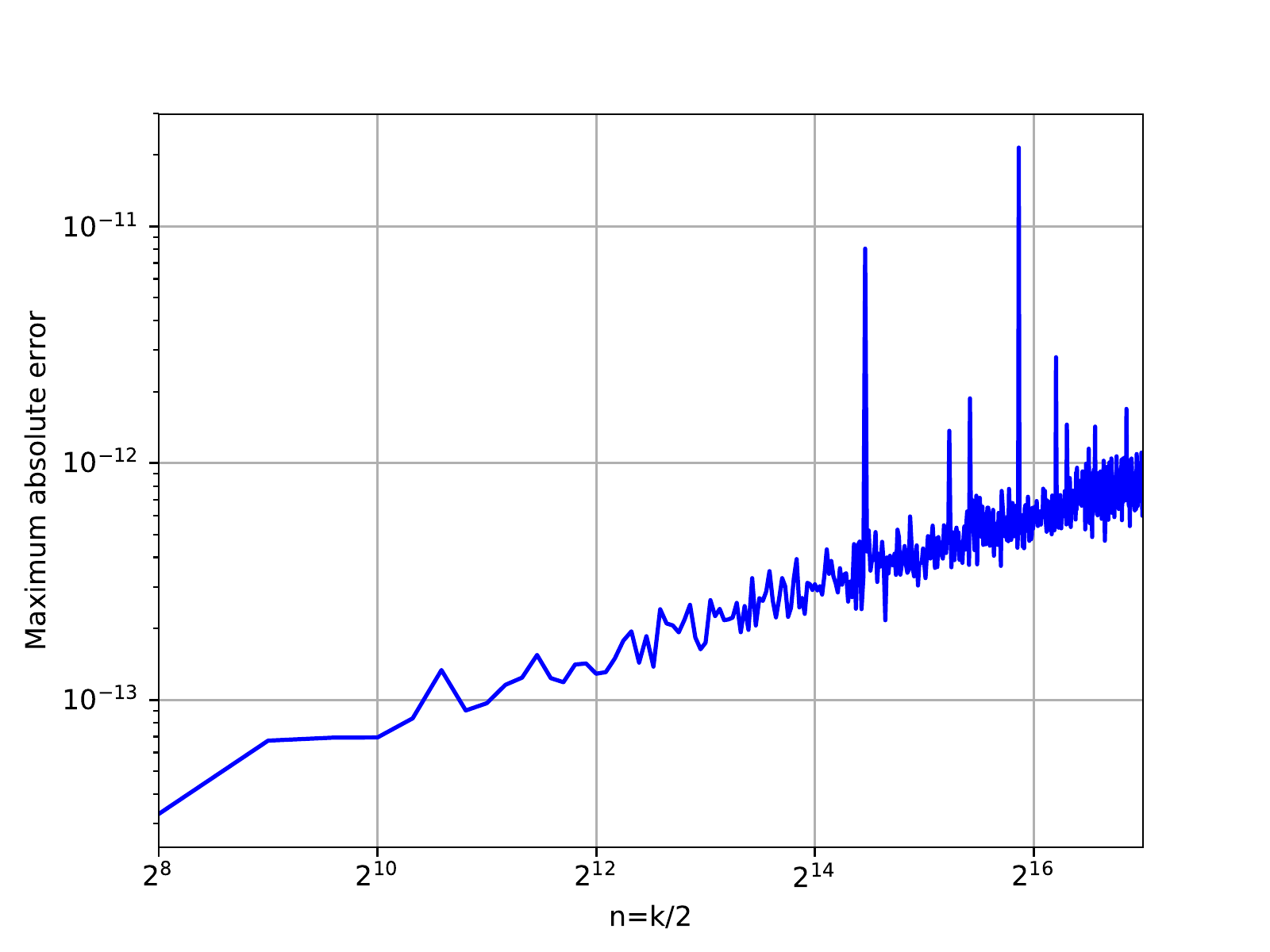}
\hfill

\hfill
\includegraphics[width=.40\textwidth]{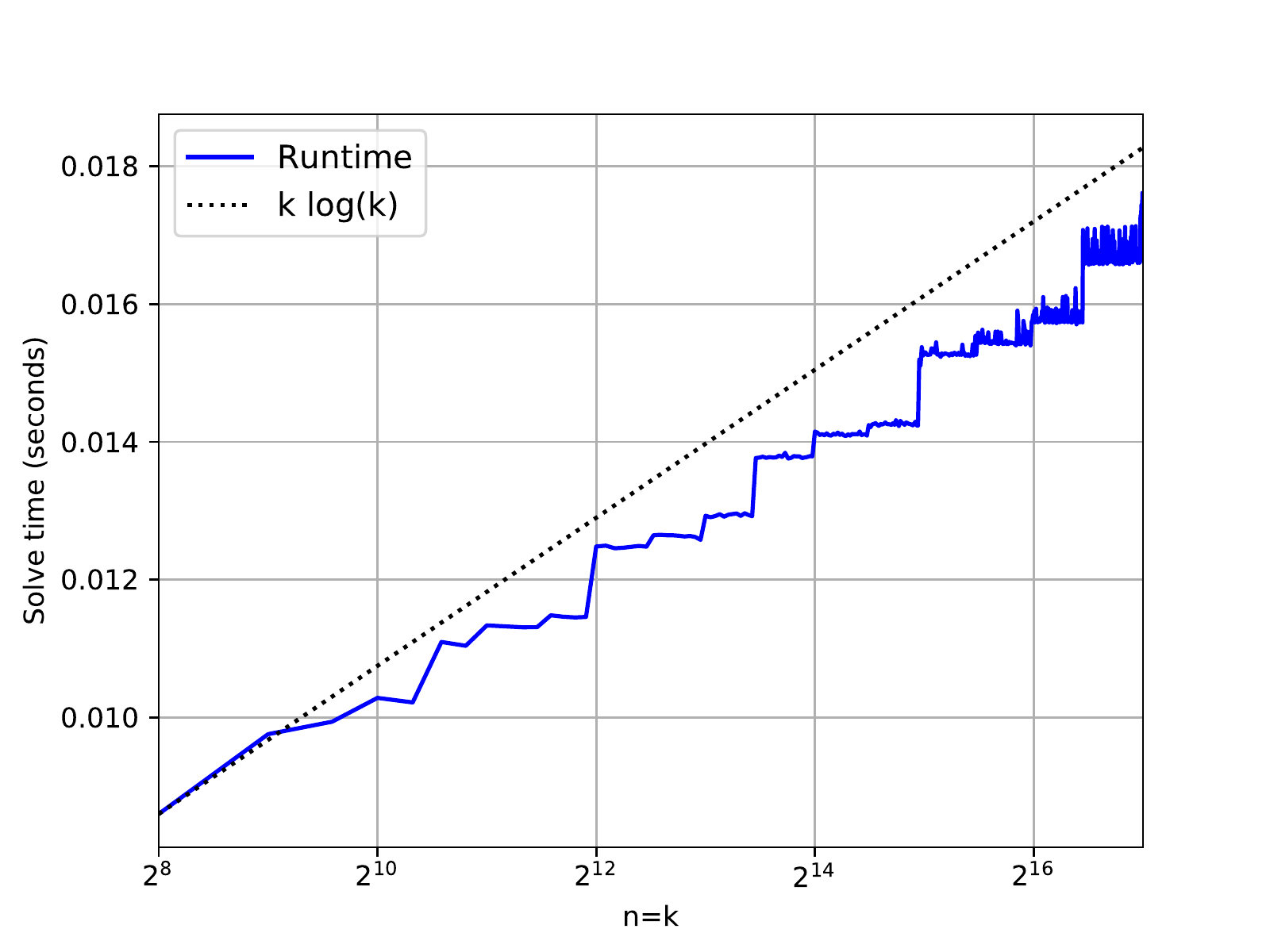}
\hfill
\includegraphics[width=.40\textwidth]{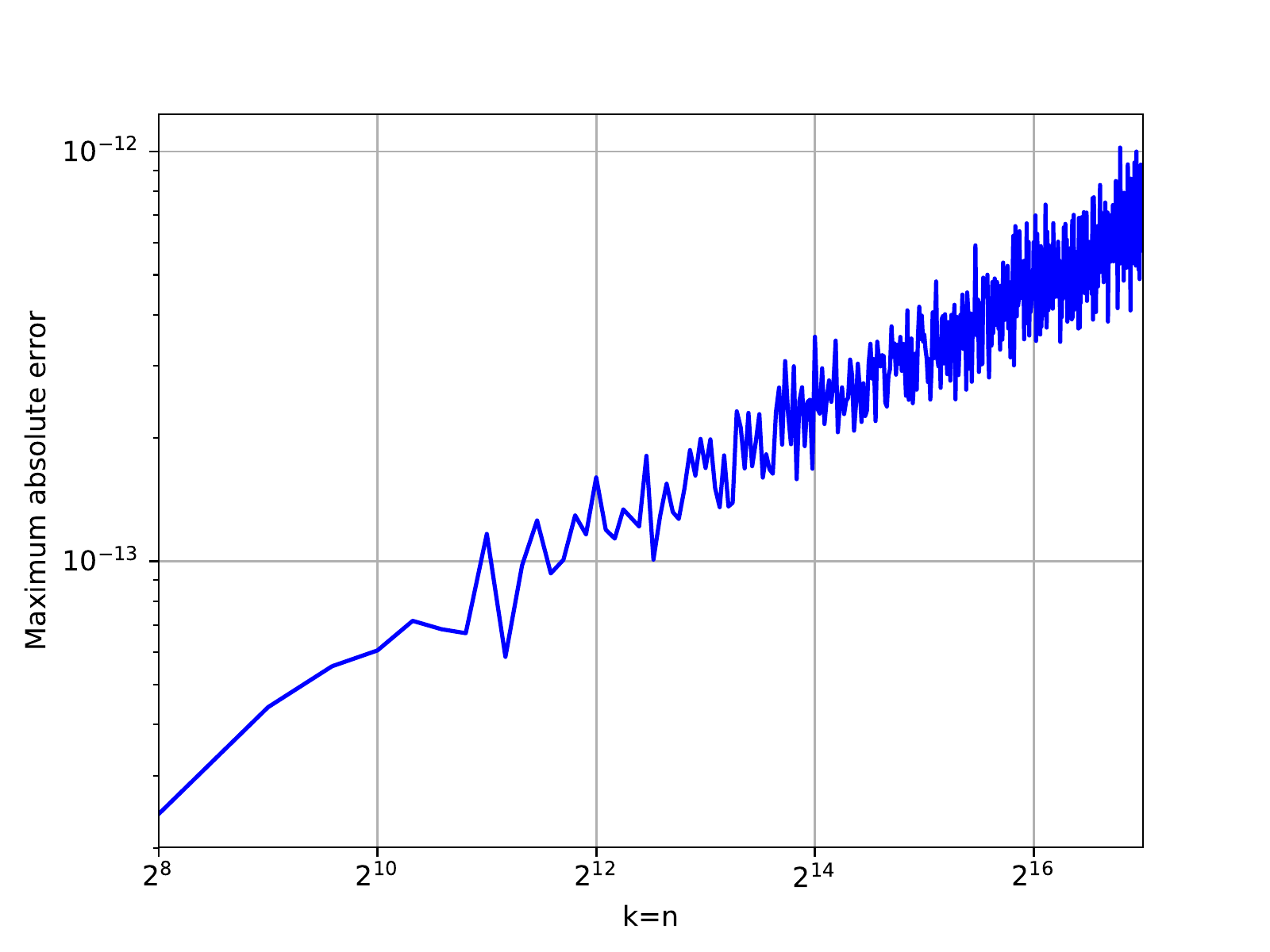}
\hfill

\caption{
The results of the experiments described in Section~\ref{section:experiments:pbessel1}
in which the performance of our solver for the perturbed Bessel equation was
tested.  In each of these experiments $q(r) = r^2-1$ and $R=2$.
The first row gives the running time and largest observed absolute error 
 as functions of $n$ when $k$ is held fixed at $2^{17}$.
The second row gives them as functions of $k$ when  $n=0$ and $k$ is increased
from $0$ to $2^{17}$.    The third and fourth row give the running time
and largest observed absolute error as functions of $n$ when $n$ is equal to a multiple of
$k$; $n = k/2$ in third row and $n=k$ in the fourth.
}
\label{figure:pbessel1}
\end{figure}

\begin{section}{Numerical Experiments}
\label{section:experiments}

In this section, we describe numerical experiments which were conducted to evaluate
the performance of the algorithm of this paper.
Our code was written in Fortran with OpenMP extensions
and was compiled with the GNU Fortran compiler version 7.4.0.
All calculations were performed on a workstation computer equipped with
$28$ Intel Xeon E5-2697 processor cores running at 2.6 GHz.
We used  P. Swarztrauber's {FFTPACK} library \cite{fftpack}  to apply the fast Fourier transform.  
We used a code  provided by V. Rokhlin to evaluate the Hankel functions.
Our implementation of the algorithm of this paper and 
our code for conducting the numerical  experiments described here 
are available on GitHub at the following address:
\begin{center}
\url{https://github.com/JamesCBremerJr/HelmRad}
\end{center}


We used the following procedure to measure the accuracy of solutions
produced by our solver for (\ref{introduction:scatbvp}).     We first  executed it
using extended precision (Fortran REAL*16) arithmetic,  which gives about 33 decimal digits of accuracy.
When possible, we then used a spectral method to verify that the obtained
scattered field satisfies the partial differential equation
\begin{equation}
\Delta u_s(x) + k^2 \left(1+q(x)\right) u_s(x) =  -k^2 q(x) u_i(x) 
\label{experiments:spectral}
\end{equation}
to at least 15 decimal digits of accuracy.  Finally, we executed our  algorithm  a second time
using double precision (Fortran REAL*8) arithmetic and measured the
error in the obtained solution by comparison with the reference
solution produced using extended precision arithmetic.
The condition number of the spectral discretization of 
(\ref{experiments:spectral}) increases rapidly with the number of discretization nodes needed
and hence with $k$, which  is why 
extended precision arithmetic was necessary to verify the reference solutions.
Even so,  ill-conditioning  limited the use of this
technique to problems in which $k$ was less than or equal to $512$.  When reporting
errors, we use parentheses to indicate experiments in which we could
not verify the accuracy of our extended precision solution through a spectral
method.

The code for the precomputation phase of our algorithm is multithreaded in order to take advantage 
of the embarrassingly parallel nature of the calculation
(each $\psi_n$ can constructed entirely  independent of the others).
The FFTPACK library, on the other hand, is single-threaded and  although  the solution phase of our 
algorithm could no doubt be accelerated 
by switching to a  multithreaded FFT library, 
we opted not to do so because
the source code for the FFTPACK library was readily available and easy to modify 
to use extended precision arithmetic.

In the course of conducting these experiments, we found that there is
a large jump in the cost of applying the FFT using the FFTPACK library
when the dimension of the transform is increased from $2^{15}$ to $2^{16}$,
and this is reflected in the timings for the solution phase of our algorithm.
Since algorithms for applying the fast Fourier transform are not our principal
concern here, we did not extensively investigate this issue.
However, we suspect that it is a cache effect.

%

\begin{subsection}{The numerical solution of the perturbed Bessel equation, part I}
\label{section:experiments:pbessel1}

We now describe a set of  experiments conducted to measure performance of our method for 
solving the  normal form of the perturbed Bessel equation (\ref{pbessel:ode}). 
In each of them, $q$ was taken to be
\begin{equation}
q(r) = r^2-1,
\end{equation}
and  the equation was solved over the interval $[0,2]$.
   The functions
\begin{equation}
\left\{J_{\frac{n}{2}}\left(\frac{\lambda}{2} r^2 \right) \sqrt{r},\
Y_{\frac{n}{2}}\left(\frac{\lambda}{2} r^2 \right) \sqrt{r}\right\}
\end{equation}
form a basis in the space of solutions of (\ref{pbessel:ode}) in this case,
so any solution which is regular at the origin is necessarily multiple of
$J_{\frac{n}{2}}\left(\frac{\lambda}{2} r^2 \right) \sqrt{r}$.
This made assessing the accuracy of obtained solutions possible --- we
did so by measuring their absolute error at $100$ points on the interval $[0,2]$.

In the first of these experiments, the results of which are shown in the
first row of Figure~\ref{figure:pbessel1}, we held $k$ fixed at $2^{17} = 131,072$
and increased $n$ from $0$ to $2^{17}$.  We report the time taken
by our solver and the largest observed absolute error as functions 
of  $n$.  

In the second, the results of which are shown in the second row
of Figure~\ref{figure:pbessel1}, we let $n=0$ and increased $k$ from $0$ to $2^{17}$.
We once again report the time taken by our solver and the 
largest observed absolute error as functions  of  $k$.

In the third and fourth experiments, the results of which are reported in the third
and fourth rows of Figure~\ref{figure:pbessel1}, respectively, we 
fixed $n$ to be a constant multiple of $k$ and increased $k$ from $0$ to
$2^{17}$.  We again measured the running time and largest absolute error.
In the third experiment, we set $n=k/2$ and in the fourth we set $n=k$.

We observe that in the case in which the solutions are purely oscillatory (i.e., when $n=0$),
the running time of the procedure is essentially independent of $k$.
This is consistent with the estimates of \cite{Bremer-Rokhlin}.
When $k$ is fixed and $n$ is increased,   there is modest growth in
the runtime of the procedure.  In the cases in which the equation has a turning point
and $n$ is increased in proportion to $k$, 
 the running time of the procedure appears to grow  logarithmically with $k$.
This is consistent with our conjecture that the running time of this procedure
grows logarithmically with $k$ in the worst case.

\end{subsection}

\begin{subsection}{The numerical solution of the perturbed Bessel equation, part II}
\label{section:experiments:pbessel2}

  In most applications, it is necessary
to calculate the solutions of the perturbed Bessel equation for a fixed $k$ and
a range of value of $n$.  For instance, it is often necessary to construct
the set
\begin{equation}
S_k =\left\{ \varphi_n : n=0,1,\ldots,k \right\},
\end{equation}
where, for each $n$, $\varphi_n$ is a solution of (\ref{pbessel:ode}).
We conducted several experiments to measure the time required by our
solver to construct the set $S_k$ as a function of $k$.
In each of them, $R$ was taken to be $2$ and,
for each $k=2^8, 2^9, \ldots, 2^{16}, 2^{17}$, we measured the time
required to construct $S_k$.    In the first experiment, the results of which are shown
in Figure~\ref{figure:pbessel2:1}, $q(r)$ was taken to be $r^2-1$.  
In the second experiment, the results of which appear in Figure~\ref{figure:pbessel2:2},
$q(r) = 14 r^2 \exp(-5 r^2)$.  In the third experiment, $q(r) = 3 \chi_{[1,2]}(r)$, where
$\chi_{[1,2]}(r)$ denotes the characteristic function of the interval $[1,2]$.
The results of this third experiment appear in Figure~\ref{figure:pbessel2:3}.
In the case of the first and third experiment, the solution of (\ref{pbessel:ode})
is known, and we were able to measure the absolute errors in the obtained solutions.

We observe first that the time required to construct $S_k$ appears to be
 $\mathcal{O}\left(k \log(k)\right)$, which is consistent with our conjecture.
We also observe that the accuracy of the obtained solutions deteriorates 
as $k$ increases.  This is expected as the condition number 
of both (\ref{introduction:scatbvp}) and (\ref{pbessel:ode})
increase with $k$.

\end{subsection}

\begin{subsection}{Scattering of a plane wave from a Gaussian potential}

In this experiment, the incident field was  the plane wave
\begin{equation}
u_i(r,t) = \exp\left(i k r \cos \left( t - \frac{\pi}{4} \right)\right)
\end{equation}
and the scattering potential was defined by
\begin{equation}
q(r) = \exp(-5 r^2).
\end{equation}
We solved (\ref{introduction:scatbvp}) for each $k=2^{4}, 2^5, \ldots, 2^{17}$
and measured the time required by each phase of our algorithm.
We also measured the absolute error in the obtained solution in the fashion
described at the beginning of this section.  The results
are shown in Figure~\ref{figure:bump1} and Table~\ref{table:bump1}.    
Figure~\ref{figure:bump1} also displays images of the incoming field, scattered
field and total field when $k=16$, as well as a plot of the function $q(r)$.
A plot of the  scattering potential as a function of the two spatial variables 
$x$ and $y$ appears in  Figure~\ref{figure:bump2}.

We observe that the running time of our solver appears to grow in line
with our conjecture --- that is, as 
$\mathcal{O}\left(k \log(k)\right)$ --- and that while the accuracy of our solver 
deteriorates with $k$, 
at least  $6$-$7$ digit accuracy is achieved all cases, even
for a problem of more than $100,000$ wavelengths in size.


\label{section:experiments:bump}
\end{subsection}

\begin{subsection}{Scattering of a circular wave from a potential resembling a volcano}

In this experiment, the incident field was  the circular wave
\begin{equation}
u_i(z) = H_0\left( k \left|z - z_0\right|\right),
\end{equation}
where $z_0 = 6i$, and the scattering potential was defined by 
\begin{equation}
q(r) = 14  r^2 \exp(-5 r^2).
\end{equation}
Again, we solved (\ref{introduction:scatbvp}) for each $k=2^{4}, 2^5, \ldots, 2^{17}$
and measured the time required by each phase of our algorithm.  
We also measured the absolute error in the obtained solution in the fashion
described at the beginning of this section.  The results
are shown in Figure~\ref{figure:volcano1} and Table~\ref{table:volcano1}.    
Figure~\ref{figure:volcano1} also displays images of the incoming field, scattered
field and total field when $k=16$, as well as a plot of the function $q(r)$.
A plot of the  scattering potential as a function of the two spatial variables 
$x$ and $y$ appears in  Figure~\ref{figure:volcano2}.

Again, we see  that the running time of our solver appears to grow
as $\mathcal{O}\left(k \log(k)\right)$, and that, in the worst case
for a problem of more than one hundred thousand wavelengths in size,
roughly  $6$-$7$ digit accuracy is obtained.

\label{section:experiments:volcano}
\end{subsection}

\begin{subsection}{Scattering of a plane  wave from a discontinuous potential}

In this experiment, the incident field was the plane wave
\begin{equation}
u_i(r,t) = \exp\left(i k r \cos \left( t - \frac{\pi}{4} \right)\right)
\end{equation}
and the scattering potential was defined by
\begin{equation}
q(r)=
\begin{cases}
1 & \ \ 0 < r <1 \\
 2 & \ \ 2 < r <3 \\
0 & \mbox{otherwise}.
\end{cases}
\end{equation}
We once again solved (\ref{introduction:scatbvp}) for each $k=2^{4}, 2^5, \ldots, 2^{17}$
and measured the time required by each phase of our algorithm.  
We also measured the absolute error in the obtained solution in the fashion
described at the beginning of this section.  The results
are shown in Figure~\ref{figure:discont1} and Table~\ref{table:discont1}.    
Figure~\ref{figure:volcano1} also displays images of the incoming field, scattered
field and total field when $k=16$, as well as a plot of the function $q(r)$.
A plot of the  scattering potential as a function of the two spatial variables 
$x$ and $y$ appears in  Figure~\ref{figure:discont2}.

Again, we see  that the running time of our solver appears to grow
as $\mathcal{O}\left(k \log(k)\right)$, and that, in the worst case
for a problem of more than one hundred thousand wavelengths in size,
roughly  $6$-$7$ digit accuracy is obtained.

\label{section:experiments:discont}
\end{subsection}

\end{section}

\begin{section}{Conclusions and Future Work}

We have developed a fast method for the numerical solution of the two-dimensional variable
coefficient Helmholtz equation in the radially symmetric case.  It is based on separation
of variables and exploits the fact  that a large class of second order differential equations
admit nonoscillatory phase functions.   Using our scheme,
problems of up to several hundred thousands wavelengths in size can be solved
with relatively high accuracy (at least 6-7 digits) in a matter of minutes
on a typical workstation computer.

Even so, there are many inefficiencies in our approach.
For instance,  a reduction in the asymptotic
running time of our method for solving the perturbed Bessel differential equation
could be achieved
through the use of standard asymptotic methods for equations of the form (\ref{introduction:ode}).
We could represent the functions $\psi_n$ near turning points
via expansions in Airy functions (see, for instance, Chapter~11 of \cite{Olver} or Chapter~4 of \cite{Fedoryuk}).
There are some numerical difficulties involved in computing the coefficients in such expansions; however,
assuming that these can be overcome, the asymptotic complexity of the precomputation
phase could be reduced to $\mathcal{O}\left(\lambda\right)$ and a dramatic speedup in the running
time of our algorithm effected.

A more challenging problem is to extend the results of this paper  to the case of 
nonsymmetric scattering potentials.   The notion of phase function
extends easily to this case. 
Indeed, if $u(x,y) = \exp(\sigma(x,y))$ solves
\begin{equation}
\Delta u(x,y) + q(x,y) u(x,y) = 0,
\label{conclusion:helmholtz}
\end{equation}
then $\sigma$ must satisfy the  two-dimensional Riccati equation
\begin{equation}
\Delta \sigma(x,y) + \left|\nabla \sigma(x,y)\right|^2 + q(x,y) = 0.
\label{conclusion:riccati}
\end{equation}
Moreover, preliminary numerical experiments suggest that, just like its one-dimensional counterpart,
(\ref{conclusion:riccati})  admits solutions which are nonoscillatory even when $q$
is of large magnitude.  It seems likely, then, that the rapid numerical
solution of various boundary value problems for (\ref{conclusion:helmholtz})
can be carried out by constructing a collection of nonoscillatory
solutions of  (\ref{conclusion:riccati}).

Suppose, for example, that $\Omega \subset \mathbb{R}^2$ is the disk of radius $R>0$ centered at $0$.
If 
\begin{equation}
\sigma_{-m},\sigma_{-m+1},\ldots,\sigma_{-1},\sigma_0,\sigma_1,\ldots,\sigma_{m-1},\sigma_m
\end{equation}
are nonoscillatory solutions of (\ref{conclusion:riccati}) such that for each $n$
the restriction of $\exp(\sigma_n(x,y))$ to the boundary of $\Omega$ is equal to
the exponential function $\exp(int)$.  If the restriction of $f$
to $\partial\Omega$ admits the expansion
\begin{equation}
\sum_{n=-m}^m a_n \exp(int), 
\end{equation}
then the solution of the  Dirichlet boundary value problem
\begin{equation}
\left\{
\begin{aligned}
\Delta u(x,y) + q(x,y) u(x,y) = 0 \ \ &\mbox{in} \ \ \Omega\\
u(x,y) = f(x,y) \ \ &\mbox{on} \ \ \partial\Omega
\end{aligned}
\right.
\label{conclusion:bvp}
\end{equation}
is
\begin{equation}
\sum_{n=-m}^m a_n \exp( \sigma_n(x,y)). 
\end{equation}
Assuming that (\ref{conclusion:riccati}) can be solved in $\mathcal{O}\left(\log(k)\right)$
time,   this would provide a method for solving 
 (\ref{conclusion:bvp}) in the general case whose running time behaves
as $\mathcal{O}\left(\log(k)\right)$.
This line of inquiry is being vigorously pursued by the author and will
be reported on at a later date.

\label{section:conclusion}
\end{section}

\begin{section}{Acknowledgments}
We thank Vladimir Rokhlin for providing us with his code for evaluating the Hankel
functions and for several useful discussions.
This work was supported in part by National Science Foundation grant DMS-1418723,
and by a UC Davis Chancellor's Fellowship.
\end{section}

\vfil
\eject
\begin{section}{References}
\bibliographystyle{acm}
\bibliography{helmrad.bib}
\end{section}


\vfil \eject

\FloatBarrier


\begin{figure}[t!]
\small
\begin{minipage}[c]{.5\textwidth}
\begin{center}
\includegraphics[width=\textwidth]{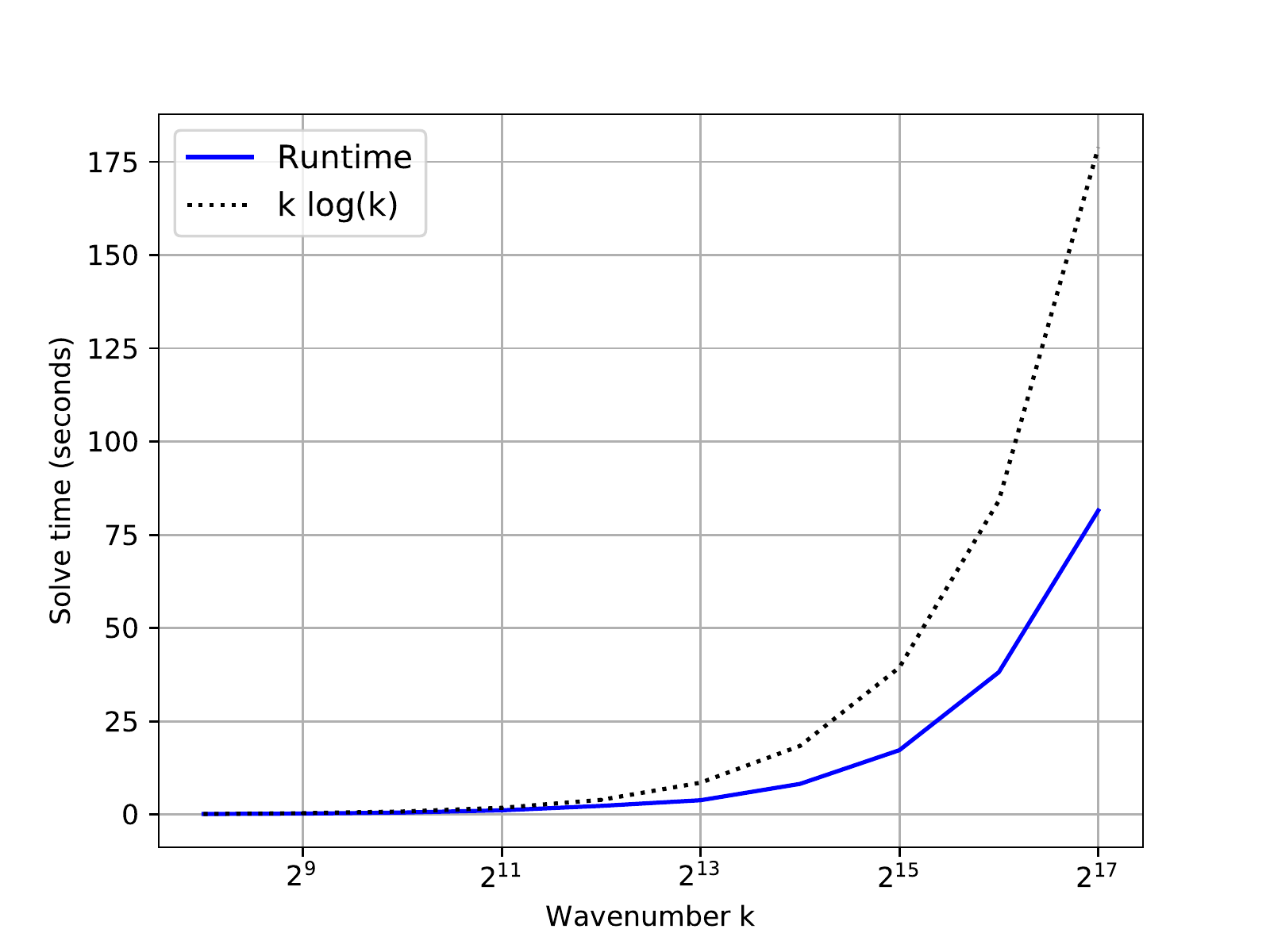}
\end{center}
\end{minipage}
\begin{minipage}[c]{.5\textwidth}
\begin{center}
\includegraphics[width=\textwidth]{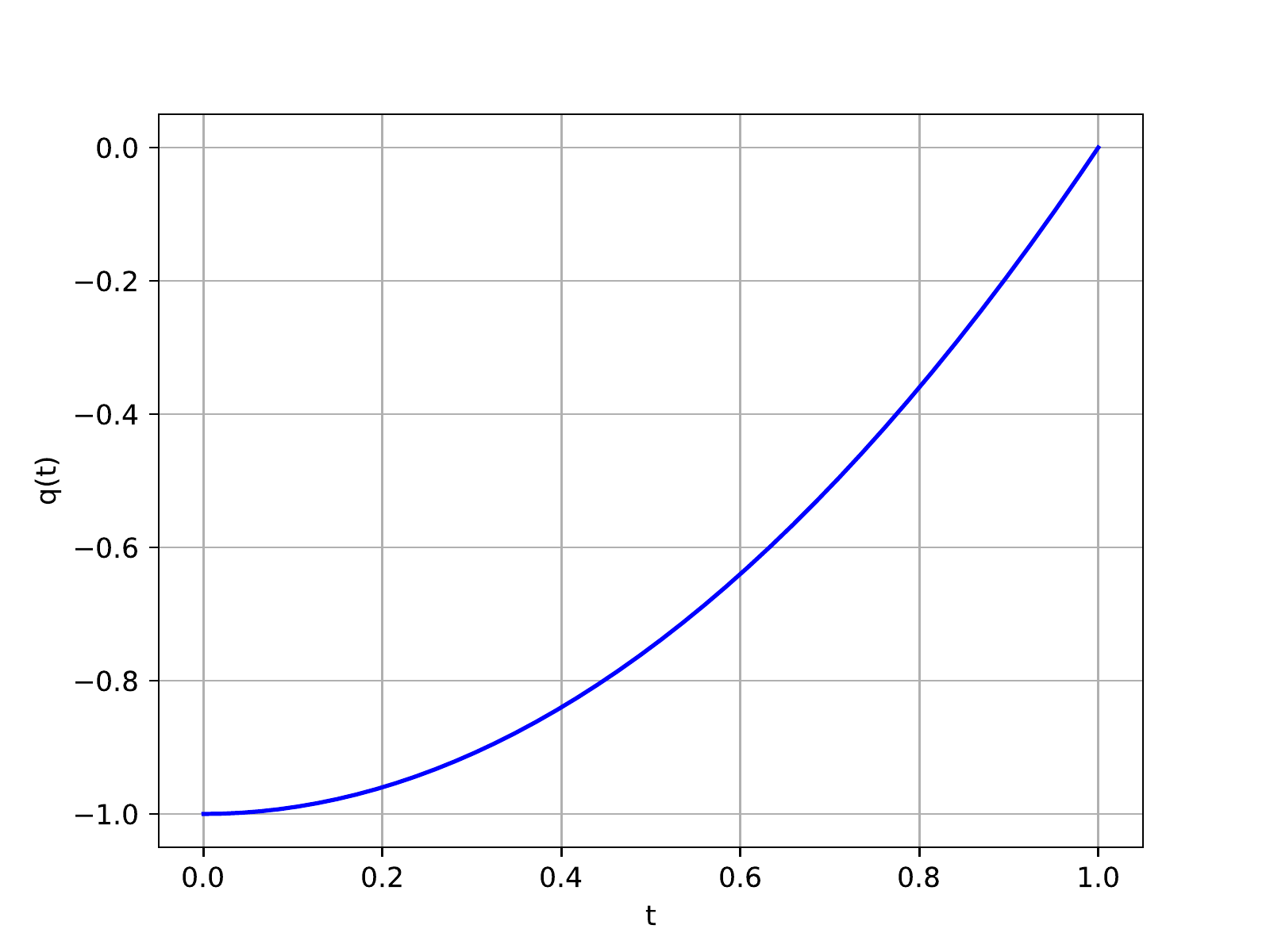}
\end{center}
\end{minipage}

\vskip 2em
\begin{center}
\begin{tabular}{cccc}
\toprule
 \addlinespace[.25em]
 \addlinespace[.125em]
$k$                                &Total time                         &Ratio of                          &Maximum absolute                          \\
                                   &(in seconds)                       &times                              &error                               \\
\midrule
 \addlinespace[.25em]
$2^{8}$ & 1.64\e{-01} &  -
 & 2.44\e{-12}  \\
 \addlinespace[.125em]
$2^{9}$ & 2.61\e{-01} & 1.59\e{+00}  & 7.69\e{-11}  \\
 \addlinespace[.125em]
$2^{10}$  & 5.38\e{-01} & 2.05\e{+00}  & 3.30\e{-12}  \\
 \addlinespace[.125em]
$2^{11}$  & 1.13\e{+00}  & 2.11\e{+00}  & 7.08\e{-12}  \\
 \addlinespace[.125em]
$2^{12}$  & 2.32\e{+00}  & 2.04\e{+00}  & 1.01\e{-11}  \\
 \addlinespace[.125em]
$2^{13}$  & 3.84\e{+00}  & 1.65\e{+00}  & 2.03\e{-11}  \\
 \addlinespace[.125em]
$2^{14}$  & 8.22\e{+00}  & 2.13\e{+00}  & 6.25\e{-10}  \\
 \addlinespace[.125em]
$2^{15}$  & 1.72\e{+01}  & 2.09\e{+00}  & 1.17\e{-10}  \\
 \addlinespace[.125em]
$2^{16}$  & 3.81\e{+01}  & 2.21\e{+00}  & 1.60\e{-09}  \\
 \addlinespace[.125em]
$2^{17}$  & 8.15\e{+01}  & 2.13\e{+00}  & 2.12\e{-09}  \\
 \addlinespace[.125em]
\bottomrule
\end{tabular}

\end{center}

\caption{
The results of the first experiment of Section~\ref{section:experiments:pbessel2}.
The plot in the upper left gives the time required to construct the set of solutions
$S_k = \{\varphi_n : n=0,\ldots,k\}$ of the perturbed Bessel equation when $q(r)  = 3 \chi_{1,2}(r)$
as function of $k$.  In the upper right is a  graph of the function $q(r)$.
Each row of the table corresponds to one value of $k$ and 
reports the time $t_k$ required to compute the set of solutions $S_k$,
the maximum observed absolute error and, when applicable, the ratio of $t_k$ to $t_{k/2}$.
}
\label{figure:pbessel2:1}
\end{figure}

\FloatBarrier

\begin{figure}[t!]
\small
\begin{minipage}[c]{.5\textwidth}
\begin{center}
\includegraphics[width=\textwidth]{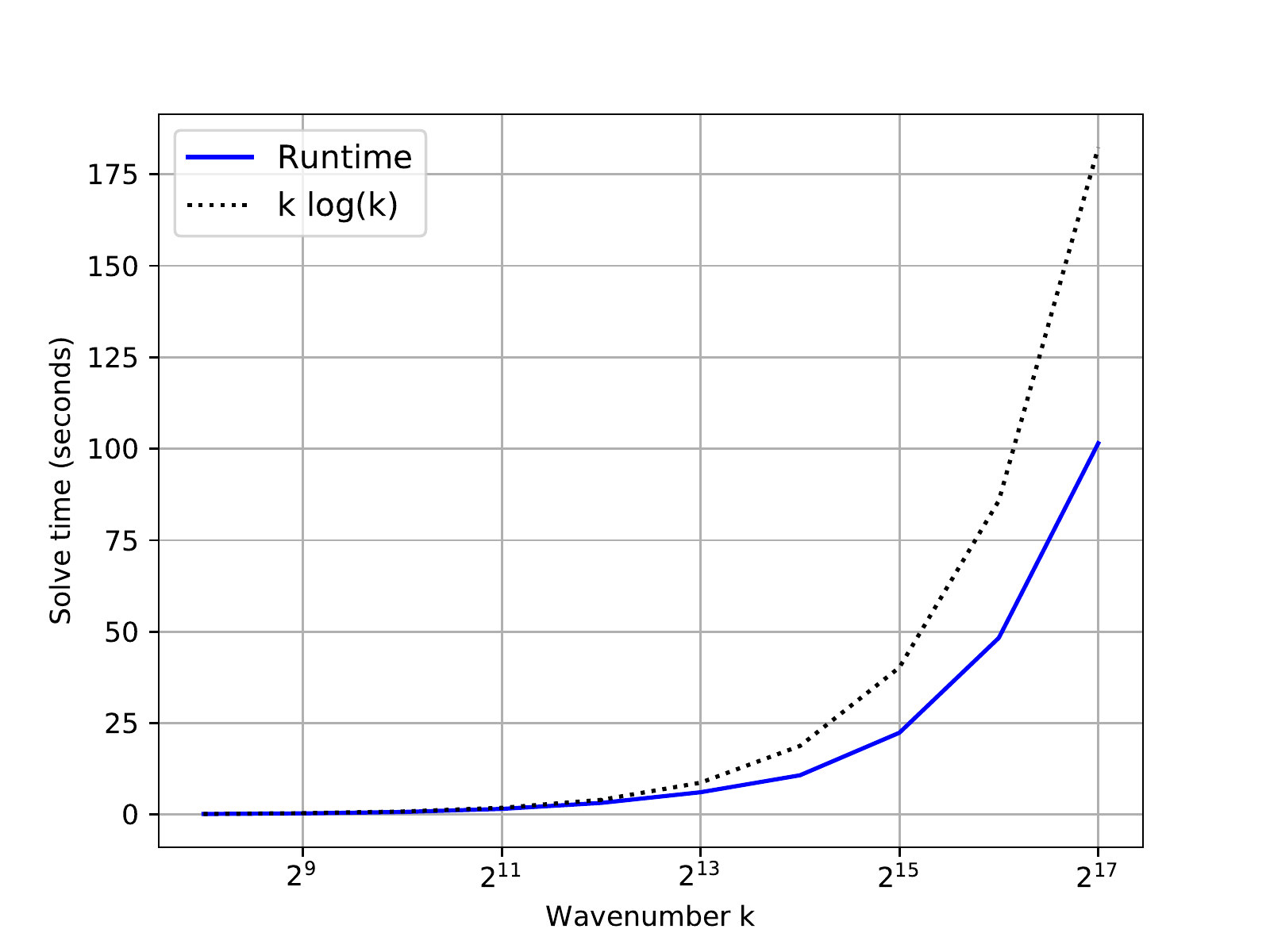}
\end{center}
\end{minipage}
\begin{minipage}[c]{.5\textwidth}
\begin{center}
\includegraphics[width=\textwidth]{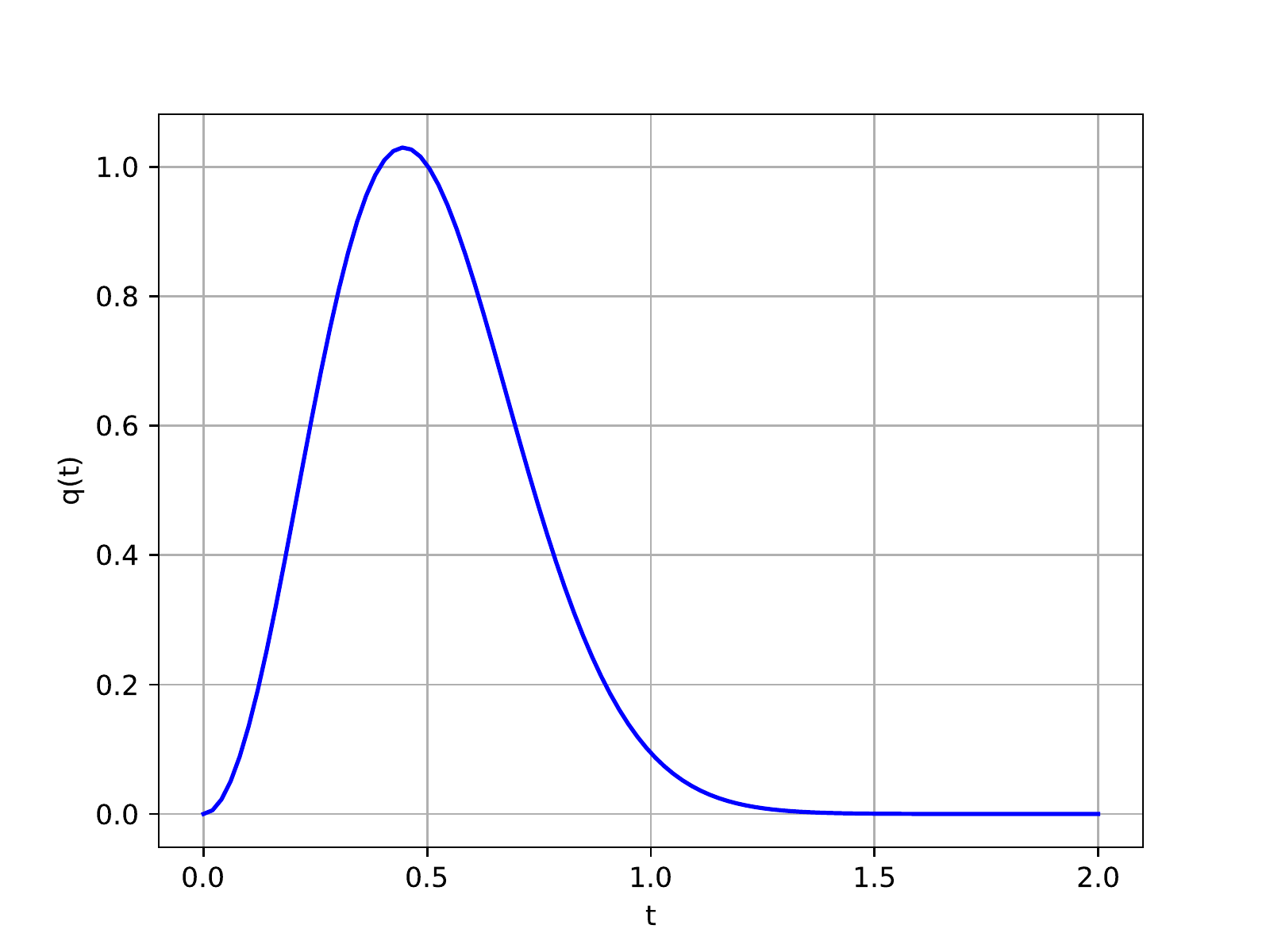}
\end{center}
\end{minipage}

\vskip 2em
\begin{center}
\begin{tabular}{ccc}
\toprule
 \addlinespace[.25em]
 \addlinespace[.125em]
$k$                                &Total time                         &Ratio of                              \\
                                   &(in seconds)                       &times                           \\
\midrule
 \addlinespace[.25em]
$2^{8}$ & 1.67\e{-01} &  -
  \\
 \addlinespace[.125em]
$2^{9}$ & 3.21\e{-01} & 1.91\e{+00}   \\
 \addlinespace[.125em]
$2^{10}$  & 7.07\e{-01} & 2.20\e{+00}   \\
 \addlinespace[.125em]
$2^{11}$  & 1.54\e{+00}  & 2.19\e{+00}   \\
 \addlinespace[.125em]
$2^{12}$  & 3.18\e{+00}  & 2.05\e{+00}   \\
 \addlinespace[.125em]
$2^{13}$  & 6.10\e{+00}  & 1.91\e{+00}   \\
 \addlinespace[.125em]
$2^{14}$  & 1.07\e{+01}  & 1.75\e{+00}   \\
 \addlinespace[.125em]
$2^{15}$  & 2.23\e{+01}  & 2.08\e{+00}   \\
 \addlinespace[.125em]
$2^{16}$  & 4.83\e{+01}  & 2.15\e{+00}   \\
 \addlinespace[.125em]
$2^{17}$  & 1.01\e{+02}  & 2.10\e{+00}   \\
 \addlinespace[.125em]
\bottomrule
\end{tabular}

\end{center}

\caption{
The results of the second experiment of Section~\ref{section:experiments:pbessel2}.
The plot in the upper left gives the time required to construct the set of solutions
$S_k = \{\varphi_n : n=0,\ldots,k\}$ of the perturbed Bessel equation when $q(r)  = 14 r^2 \exp(-5 r^2)$
as function of $k$.  In the upper right is a  graph of the function $q(r)$.
Each row of the table corresponds to one value of $k$ and 
reports the time $t_k$ required to compute the set of solutions $S_k$
as well as the ratio of $t_k$ to $t_{k/2}$ (when applicable).
}
\label{figure:pbessel2:2}

\end{figure}

\FloatBarrier


\begin{figure}[t!]
\small
\begin{minipage}[c]{.5\textwidth}
\begin{center}
\includegraphics[width=\textwidth]{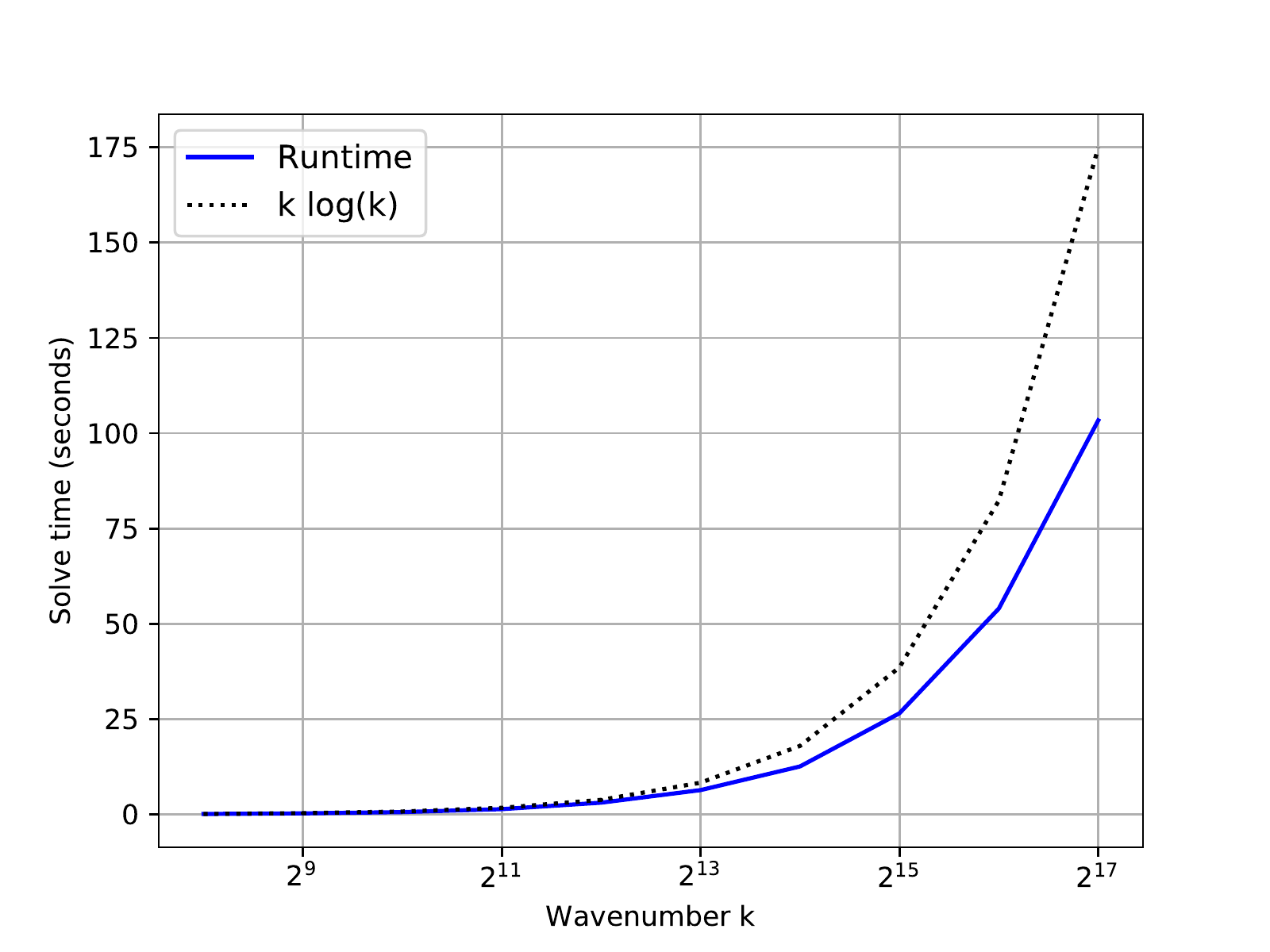}
\end{center}
\end{minipage}
\begin{minipage}[c]{.5\textwidth}
\begin{center}
\includegraphics[width=\textwidth]{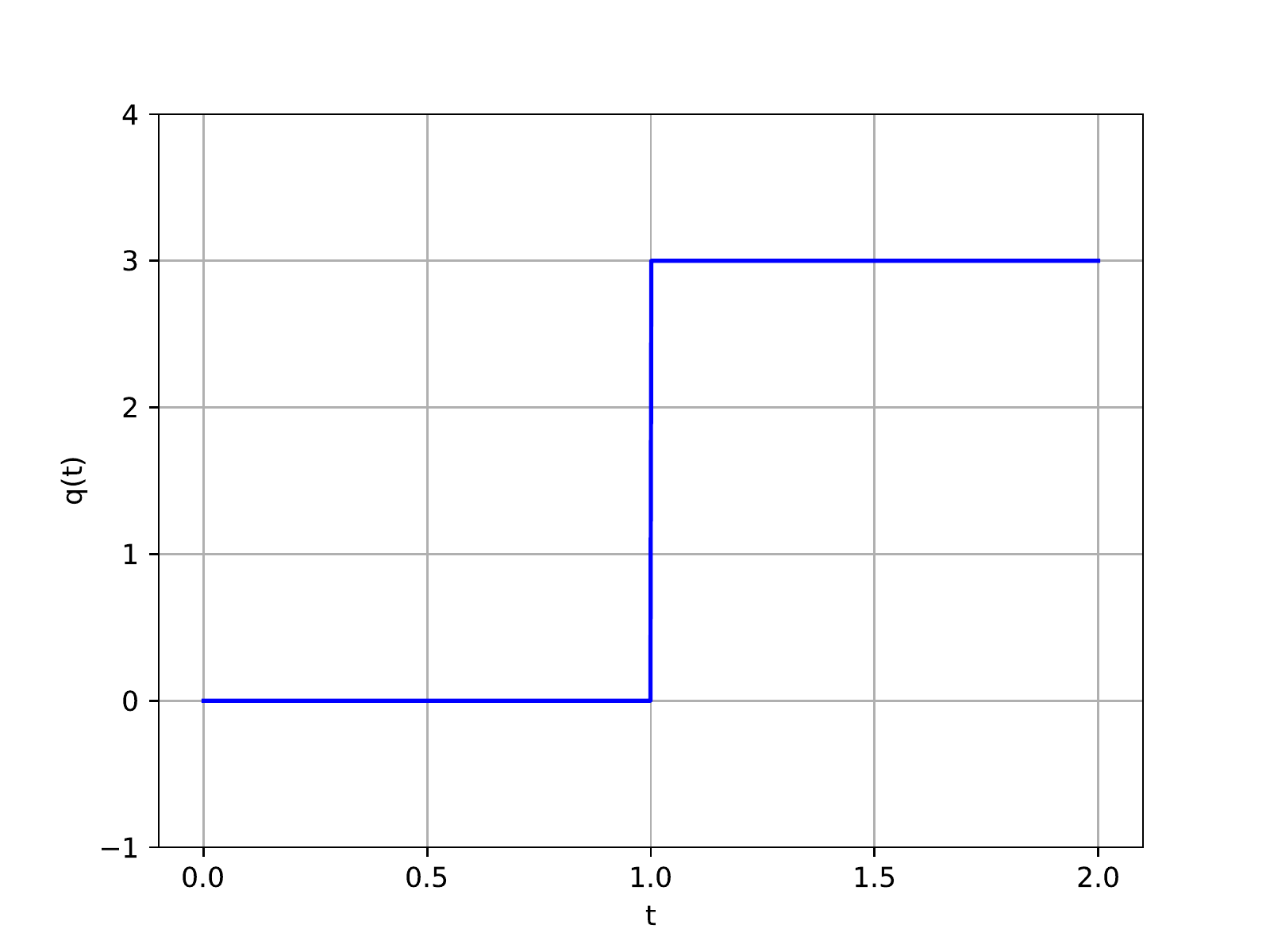}
\end{center}
\end{minipage}

\vskip 2em
\begin{center}
\begin{tabular}{cccc}
\toprule
 \addlinespace[.25em]
 \addlinespace[.125em]
$k$                                &Total time                         &Ratio of                          &Maximum absolute                          \\
                                   &(in seconds)                       &times                              &error                               \\
\midrule
 \addlinespace[.25em]
$2^{8}$ & 1.60\e{-01} &  -
 & 2.79\e{-12}  \\
 \addlinespace[.125em]
$2^{9}$ & 3.14\e{-01} & 1.95\e{+00}  & 2.23\e{-12}  \\
 \addlinespace[.125em]
$2^{10}$  & 6.50\e{-01} & 2.07\e{+00}  & 1.59\e{-12}  \\
 \addlinespace[.125em]
$2^{11}$  & 1.42\e{+00}  & 2.18\e{+00}  & 1.35\e{-12}  \\
 \addlinespace[.125em]
$2^{12}$  & 3.12\e{+00}  & 2.19\e{+00}  & 3.84\e{-12}  \\
 \addlinespace[.125em]
$2^{13}$  & 6.41\e{+00}  & 2.05\e{+00}  & 2.72\e{-11}  \\
 \addlinespace[.125em]
$2^{14}$  & 1.26\e{+01}  & 1.96\e{+00}  & 8.03\e{-11}  \\
 \addlinespace[.125em]
$2^{15}$  & 2.65\e{+01}  & 2.10\e{+00}  & 1.10\e{-10}  \\
 \addlinespace[.125em]
$2^{16}$  & 5.40\e{+01}  & 2.03\e{+00}  & 1.05\e{-09}  \\
 \addlinespace[.125em]
$2^{17}$  & 1.03\e{+02}  & 1.91\e{+00}  & 5.14\e{-09}  \\
 \addlinespace[.125em]
\bottomrule
\end{tabular}

\end{center}

\caption{
The results of the first experiment of Section~\ref{section:experiments:pbessel2}.
The plot in the upper left gives the time required to construct the set of solutions
$S_k = \{\varphi_n : n=0,\ldots,k\}$ of the perturbed Bessel equation when $q(r)  = 3 \chi_{1,2}(r)$
as function of $k$.  In the upper right is a  graph of the function $q(r)$.
Each row of the table corresponds to one value of $k$ and 
reports the time $t_k$ required to compute the set of solutions $S_k$,
the maximum observed absolute error and, when applicable, the ratio of $t_k$ to $t_{k/2}$.
}
\label{figure:pbessel2:3}
\end{figure}

\FloatBarrier


\begin{figure}[p]
\small
\begin{minipage}[c]{.5\textwidth}
\begin{center}
\includegraphics[width=\textwidth]{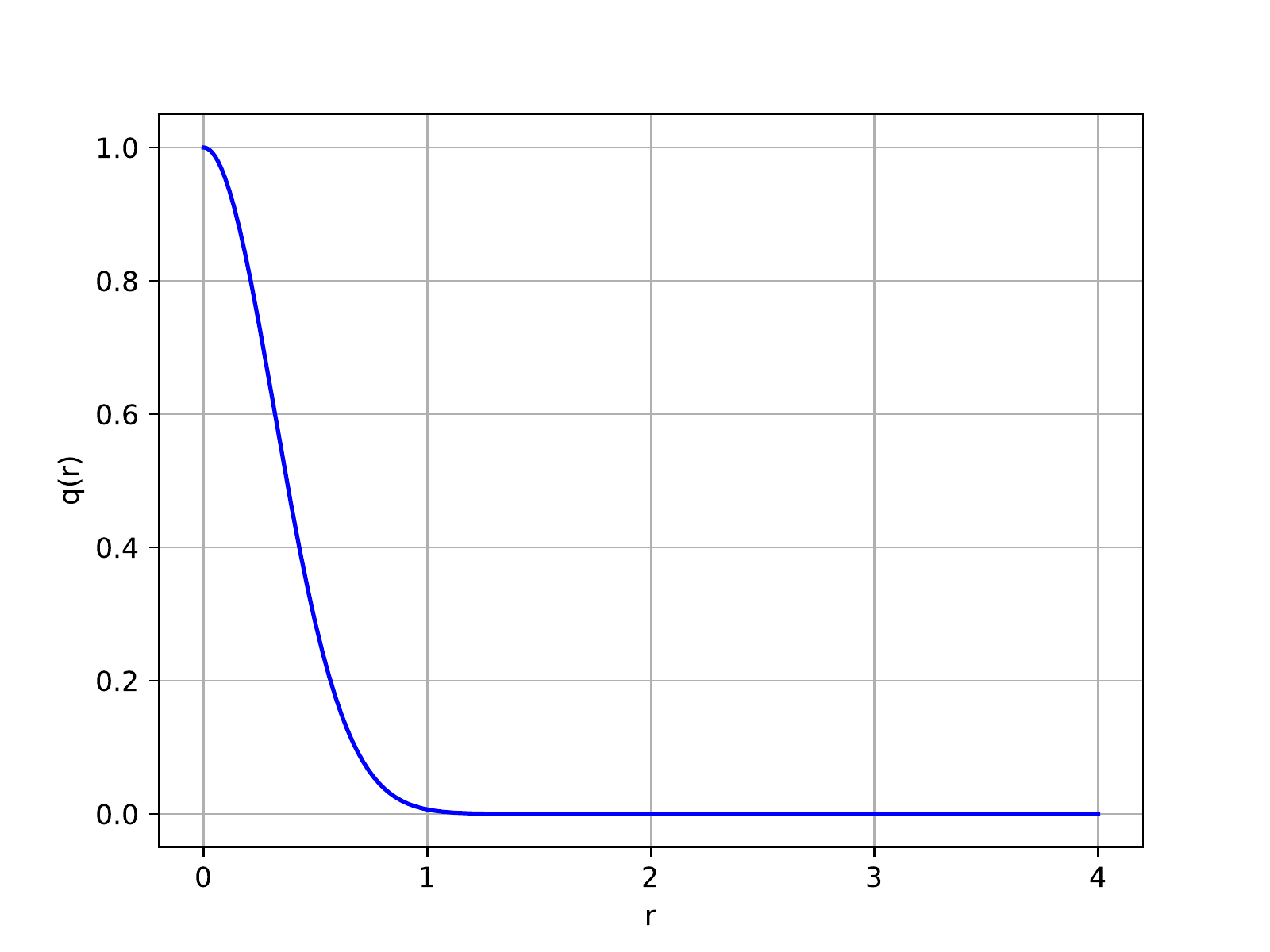}
\end{center}
\end{minipage}
\begin{minipage}[c]{.5\textwidth}
\begin{center}
\includegraphics[width=\textwidth]{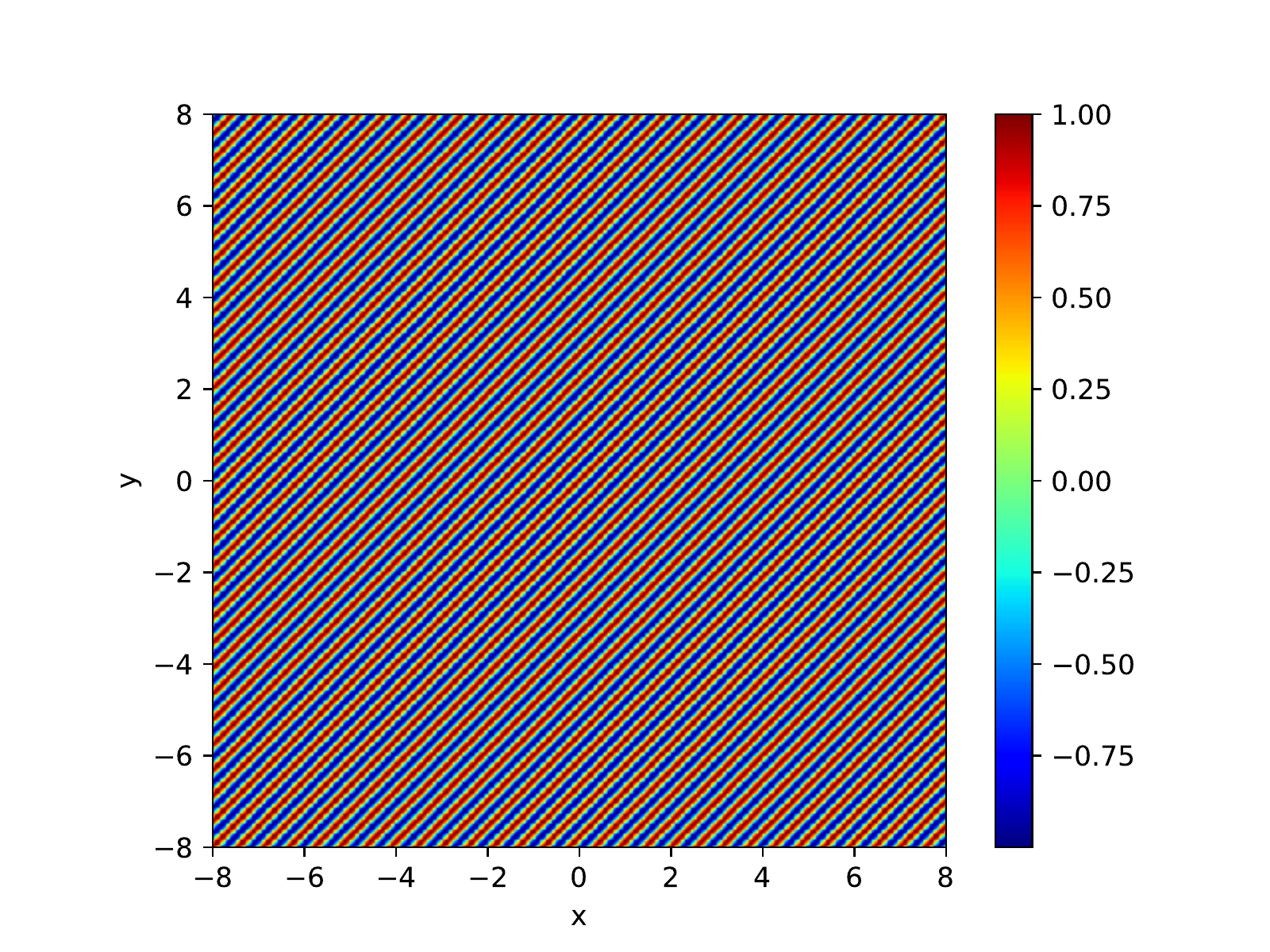}
\end{center}
\end{minipage}

\begin{minipage}[c]{.5\textwidth}
\begin{center}
\includegraphics[width=\textwidth]{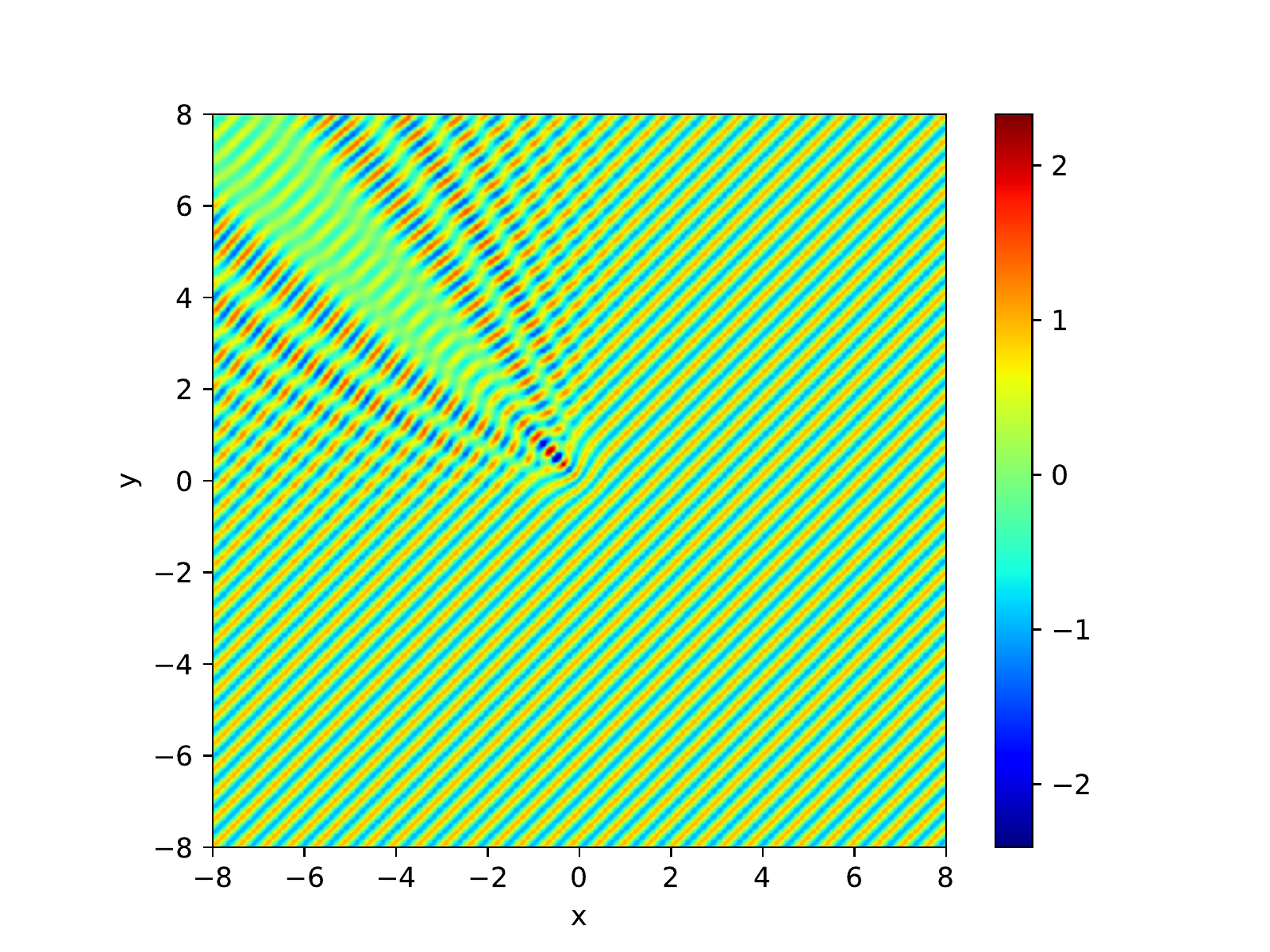}
\end{center}
\end{minipage}
\begin{minipage}[c]{.5\textwidth}
\begin{center}
\includegraphics[width=\textwidth]{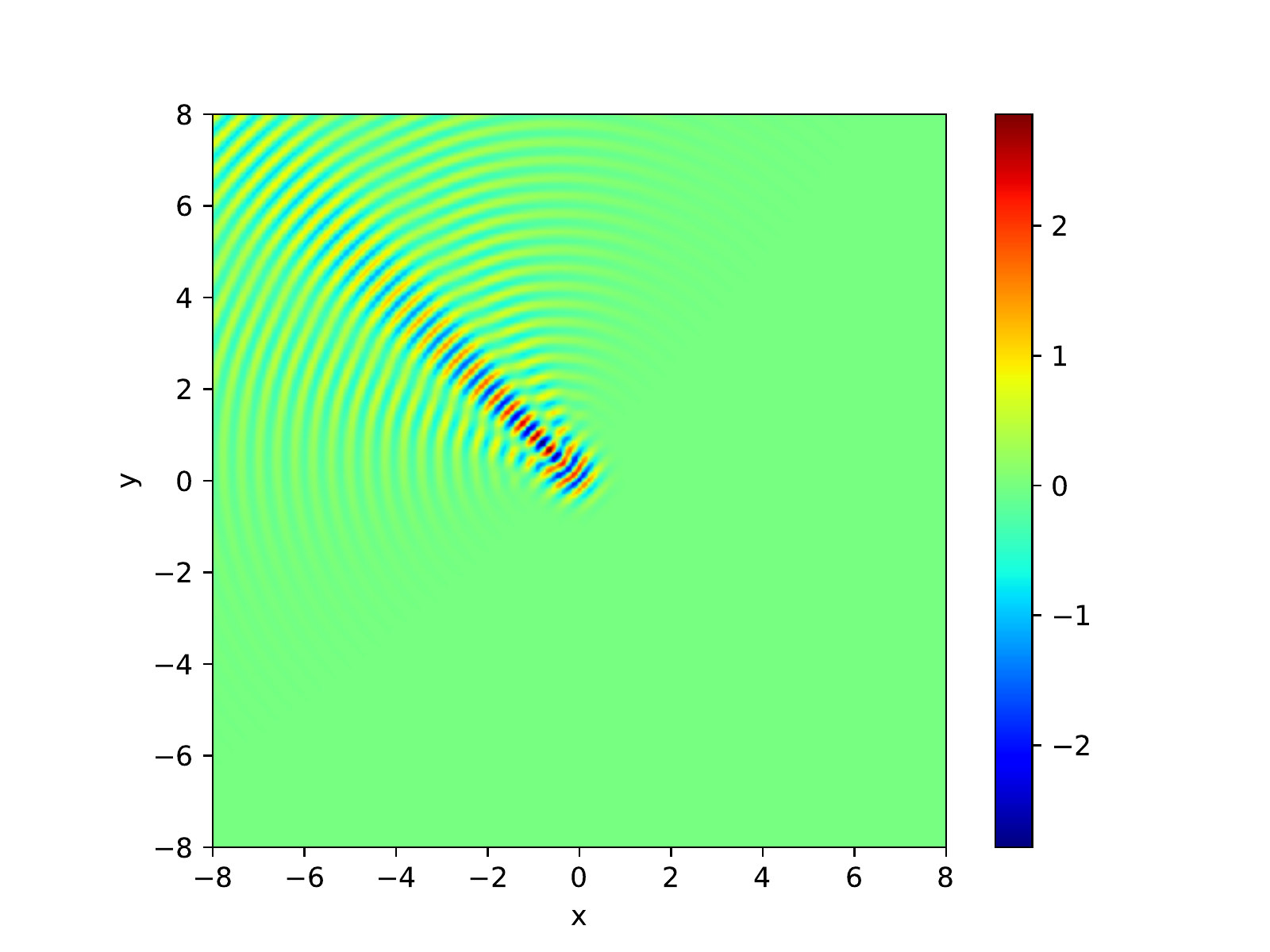}
\end{center}
\end{minipage}


\begin{minipage}[c]{.5\textwidth}
\begin{center}
\includegraphics[width=\textwidth]{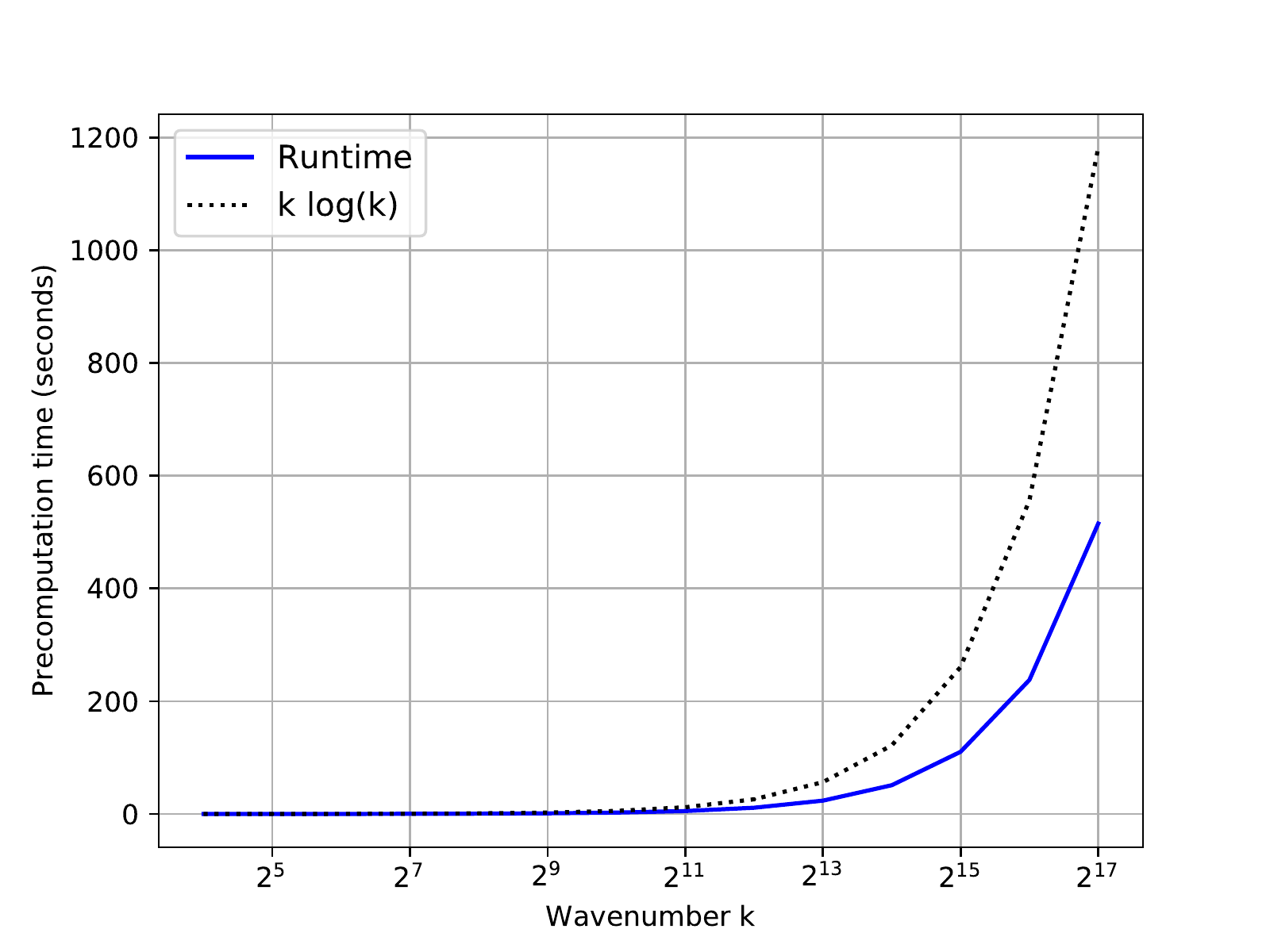}
\end{center}
\end{minipage}
\begin{minipage}[c]{.5\textwidth}
\begin{center}
\includegraphics[width=\textwidth]{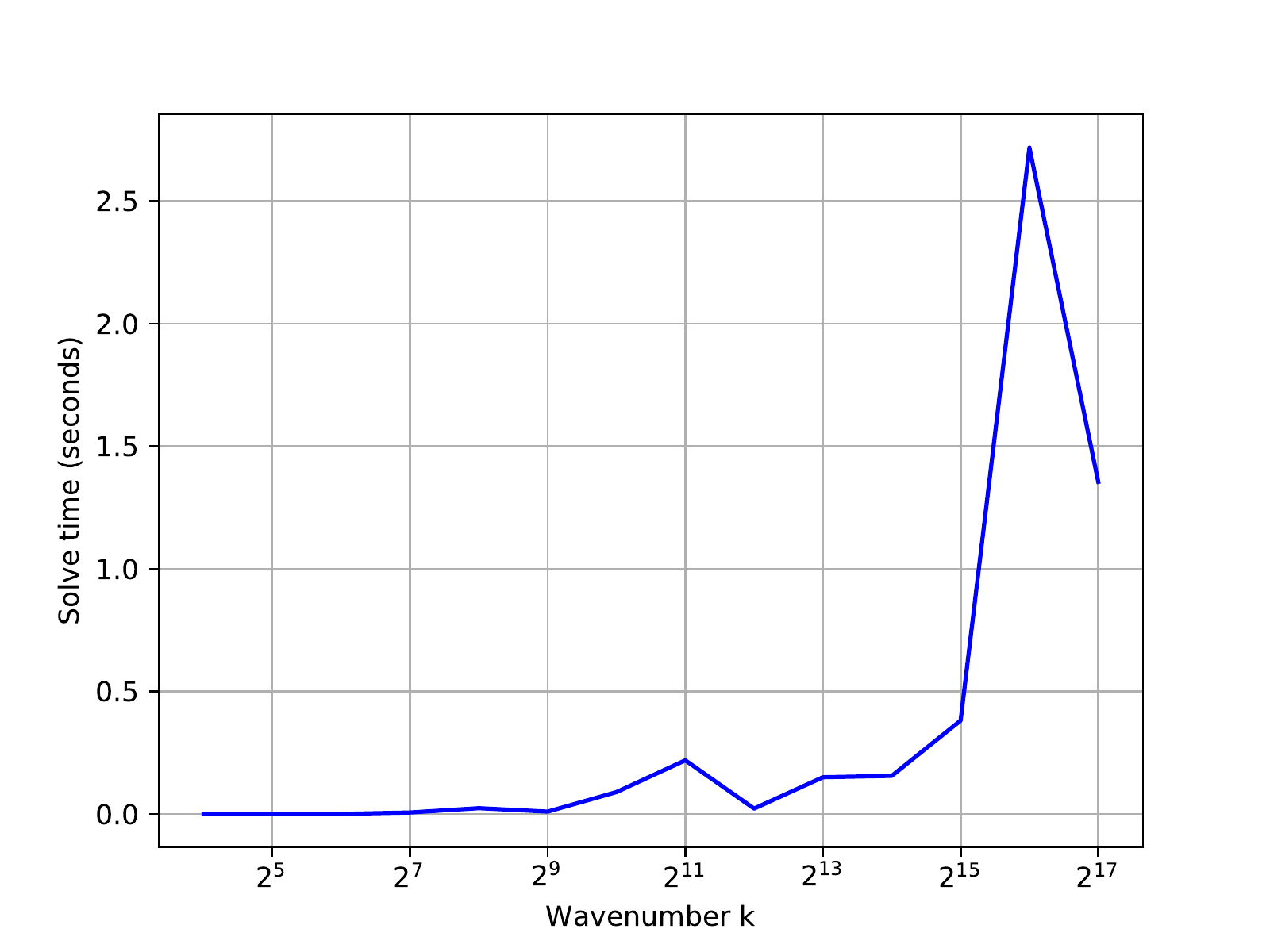}
\end{center}
\end{minipage}
\caption{Figures related to the experiment of Section~\ref{section:experiments:bump}.
In the upper left is a plot of the function $q(r)$ and in the upper right is
an image of the incident wave when $k=16$.   The image at middle left is of the
real part of the total field $k=16$.  At middle right is an image of
the real part of the scattered field when $k=16$.
At bottom left, the running time of the precomputation phase is plotted as a function of $k$
and at bottom right, the running time of the solution phase is plotted as a function of $k$.
}
\label{figure:bump1}
\end{figure}

\begin{table}[p]
\begin{center}
\begin{tabular}{crccc}
\toprule
 \addlinespace[.25em]
$k$                                &$m$                                &Maximum absolute                   &Precomp time                       &Solve time                        \\
                                   &                                   &error                              &(in seconds)                       &(in seconds)                       \\
\midrule
 \addlinespace[.25em]
$2^{4}$ &      100 & 8.34\e{-14} & 3.39\e{-02} & 1.31\e{-04}  \\
 \addlinespace[.125em]
$2^{5}$ &      201 & 1.55\e{-12} & 6.33\e{-02} & 2.34\e{-04}  \\
 \addlinespace[.125em]
$2^{6}$ &      402 & 1.26\e{-12} & 1.30\e{-01} & 2.99\e{-04}  \\
 \addlinespace[.125em]
$2^{7}$ &      804 & 1.22\e{-12} & 2.67\e{-01} & 6.15\e{-03}  \\
 \addlinespace[.125em]
$2^{8}$ &     1608 & 2.17\e{-12} & 5.51\e{-01} & 2.37\e{-02}  \\
 \addlinespace[.125em]
$2^{9}$ &     3216 & $\left(4.72\e{-12}\right)$ & 1.16\e{+00}  & 9.63\e{-03}  \\
 \addlinespace[.125em]
$2^{10}$  &     6433 & $\left(8.54\e{-12}\right)$ & 2.46\e{+00}  & 8.94\e{-02}  \\
 \addlinespace[.125em]
$2^{11}$  &    12867 & $\left(1.85\e{-11}\right)$ & 5.22\e{+00}  & 2.19\e{-01}  \\
 \addlinespace[.125em]
$2^{12}$  &    25735 & $\left(6.13\e{-11}\right)$ & 1.10\e{+01}  & 2.21\e{-02}  \\
 \addlinespace[.125em]
$2^{13}$  &    51471 & $\left(2.05\e{-10}\right)$ & 2.36\e{+01}  & 1.50\e{-01}  \\
 \addlinespace[.125em]
$2^{14}$  &   102943 & $\left(1.51\e{-09}\right)$ & 5.10\e{+01}  & 1.55\e{-01}  \\
 \addlinespace[.125em]
$2^{15}$  &   205887 & $\left(3.85\e{-09}\right)$ & 1.10\e{+02}  & 3.81\e{-01}  \\
 \addlinespace[.125em]
$2^{16}$  &   411774 & $\left(2.16\e{-08}\right)$ & 2.37\e{+02}  & 2.71\e{+00}   \\
 \addlinespace[.125em]
$2^{17}$  &   823549 & $\left(1.01\e{-07}\right)$ & 5.15\e{+02}  & 1.35\e{+00}   \\
 \addlinespace[.125em]
\bottomrule
\end{tabular}

\end{center}
\caption{The results of the experiments of Section~\ref{section:experiments:bump}. Each row
of the table corresponds to one wavenumber $k$ and gives the number $m$ of Fourier modes used
to represent the incident wave, the maximum observed absolute error in the obtained solution
(in cases in which this could be measured), and the time taken by each phase of our solver.
The absolute maximum errors are calculated via comparison with solutions generated
using extended precision arithmetic.  Parentheses are used to indicate cases in which
the wavenumber was too large for the  the extended precision
solution to be verified via a spectral method.
}
\label{table:bump1}
\end{table}

\begin{figure}[p]
\begin{center}
\includegraphics[width=.75\textwidth]{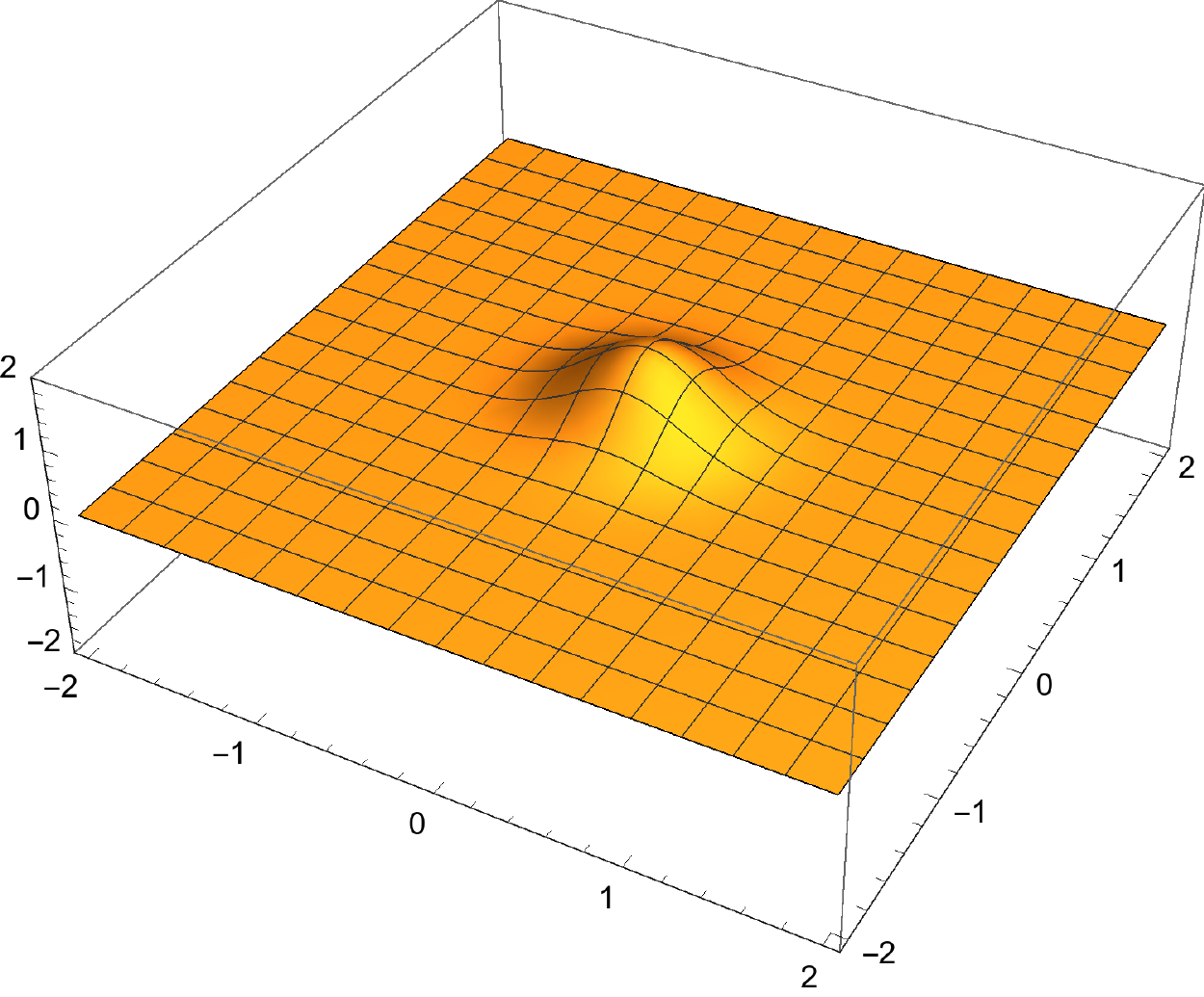}
\end{center}
\caption{A plot of the scattering potential used in the experiment of Section~\ref{section:experiments:bump}.}
\label{figure:bump2}
\end{figure}


\FloatBarrier

\begin{figure}[p]
\small
\begin{minipage}[c]{.5\textwidth}
\begin{center}
\includegraphics[width=\textwidth]{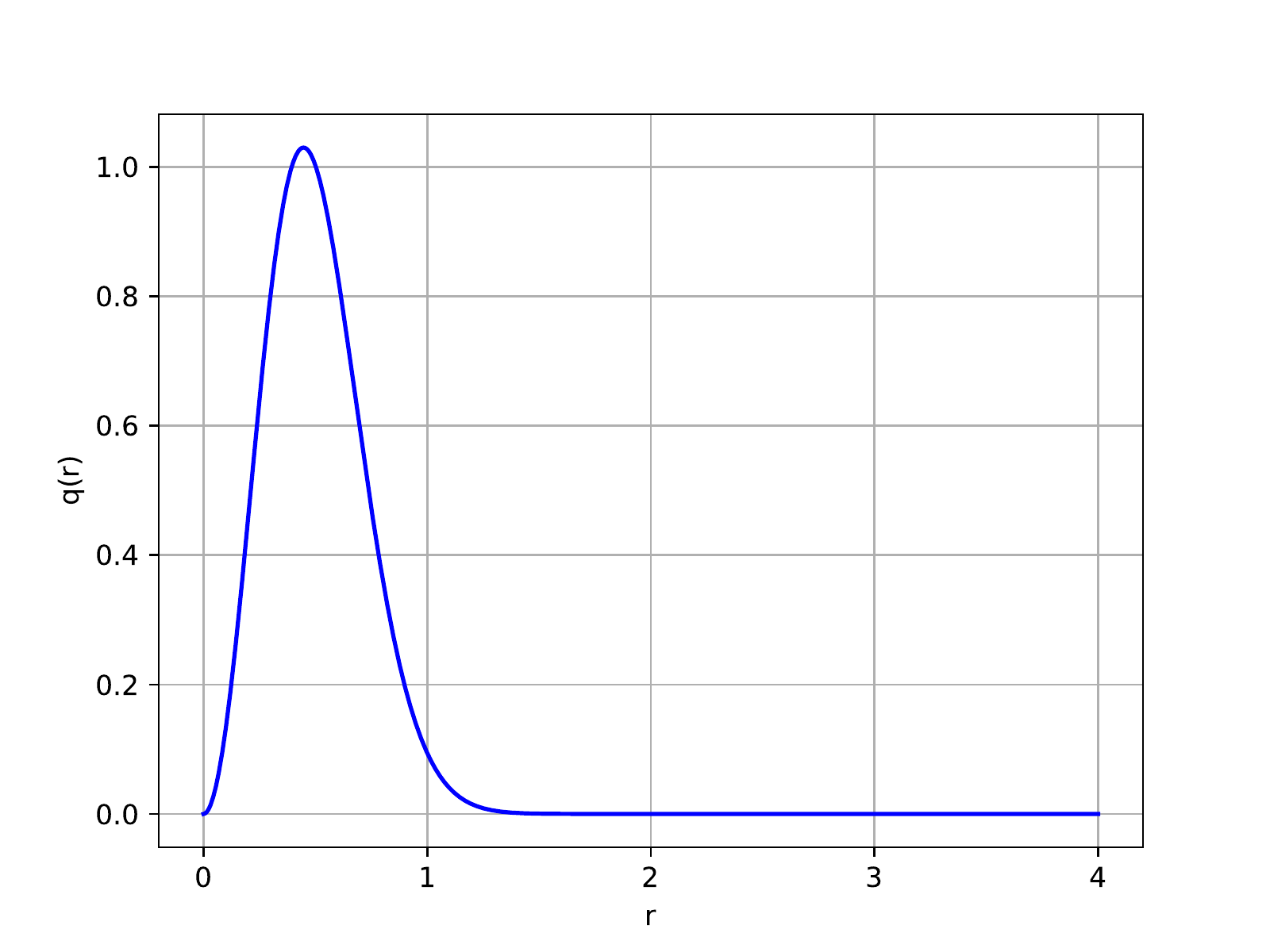}
\end{center}
\end{minipage}
\begin{minipage}[c]{.5\textwidth}
\begin{center}
\includegraphics[width=\textwidth]{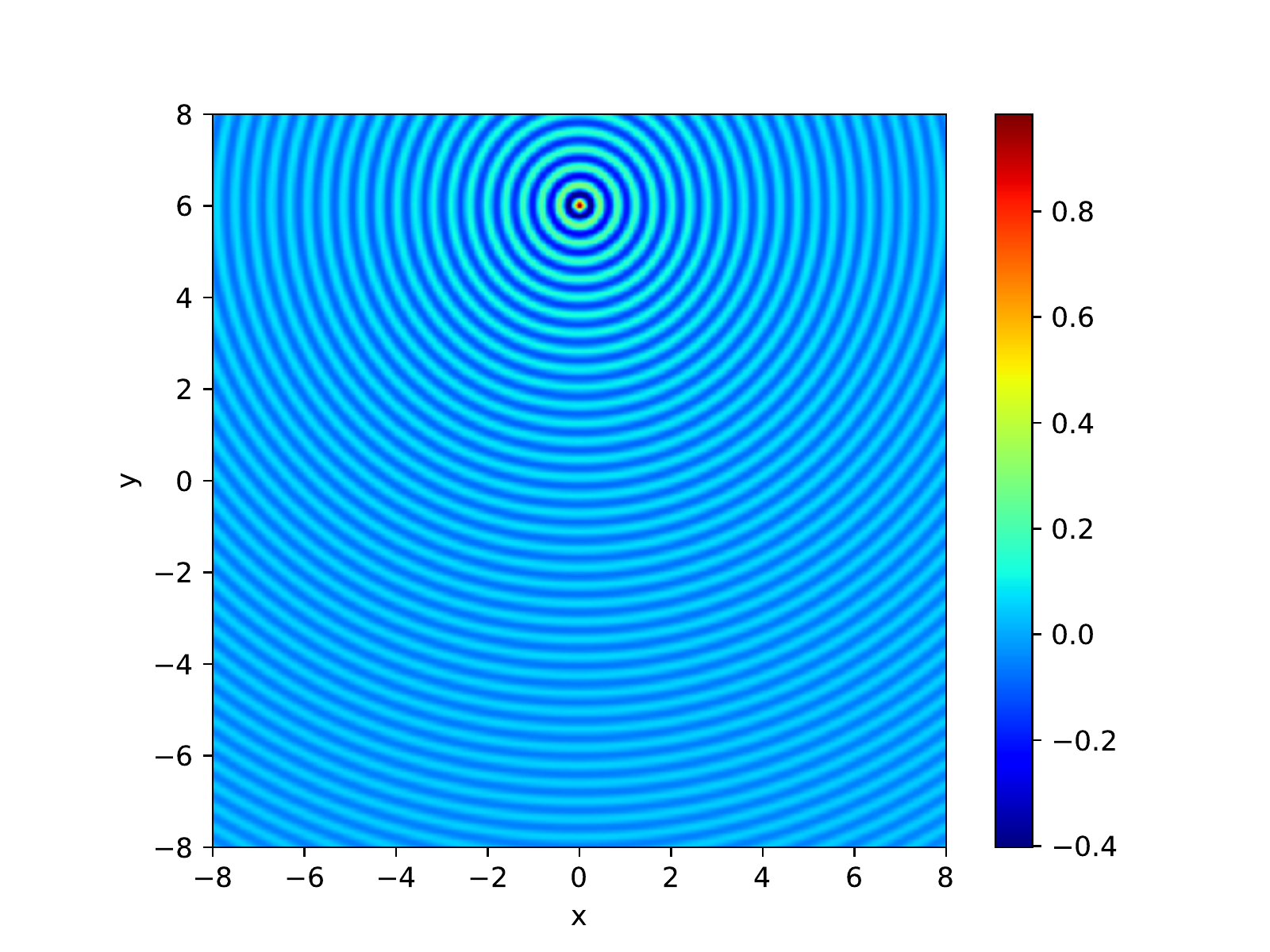}
\end{center}
\end{minipage}

\begin{minipage}[c]{.5\textwidth}
\begin{center}
\includegraphics[width=\textwidth]{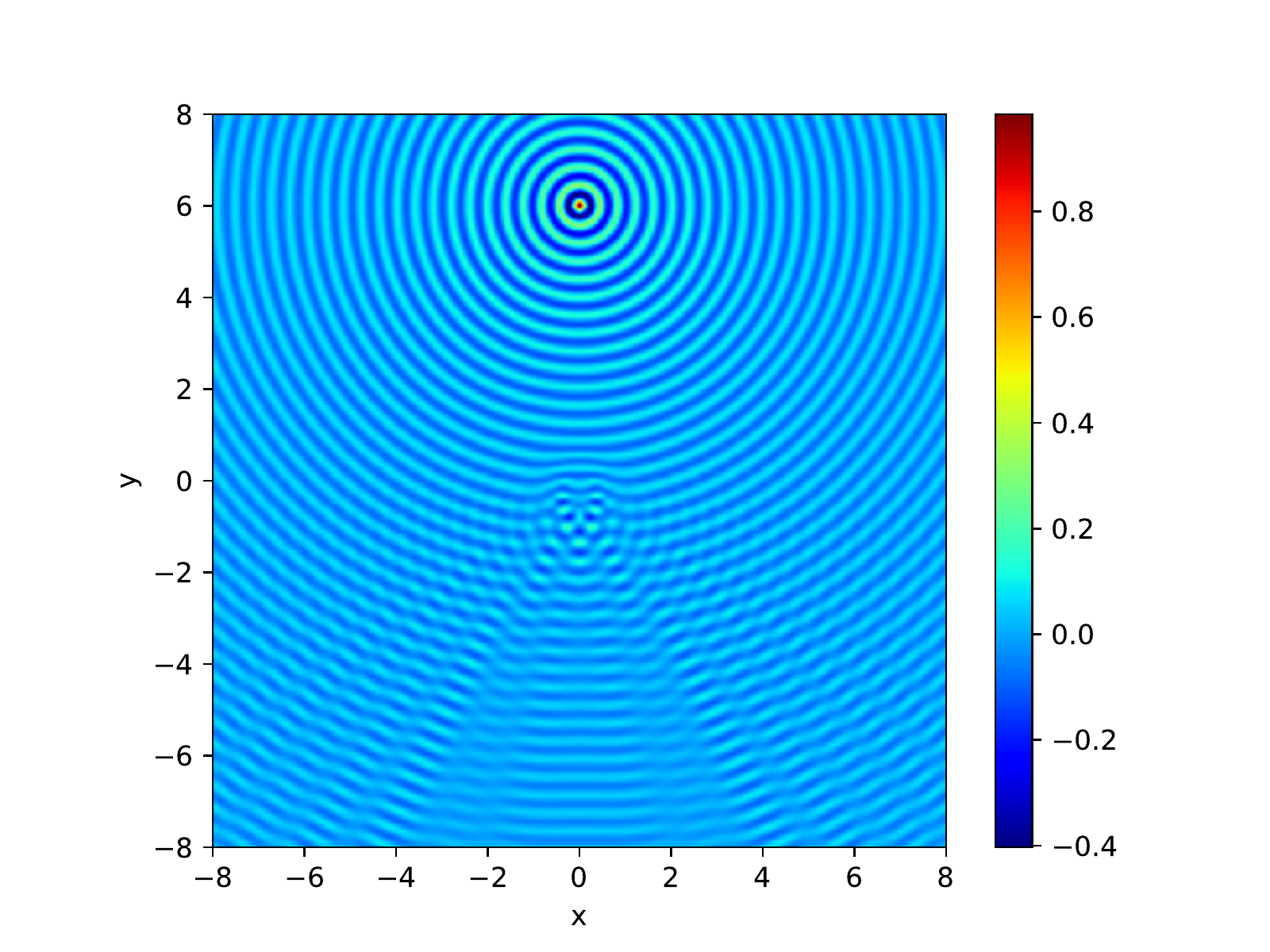}
\end{center}
\end{minipage}
\begin{minipage}[c]{.5\textwidth}
\begin{center}
\includegraphics[width=\textwidth]{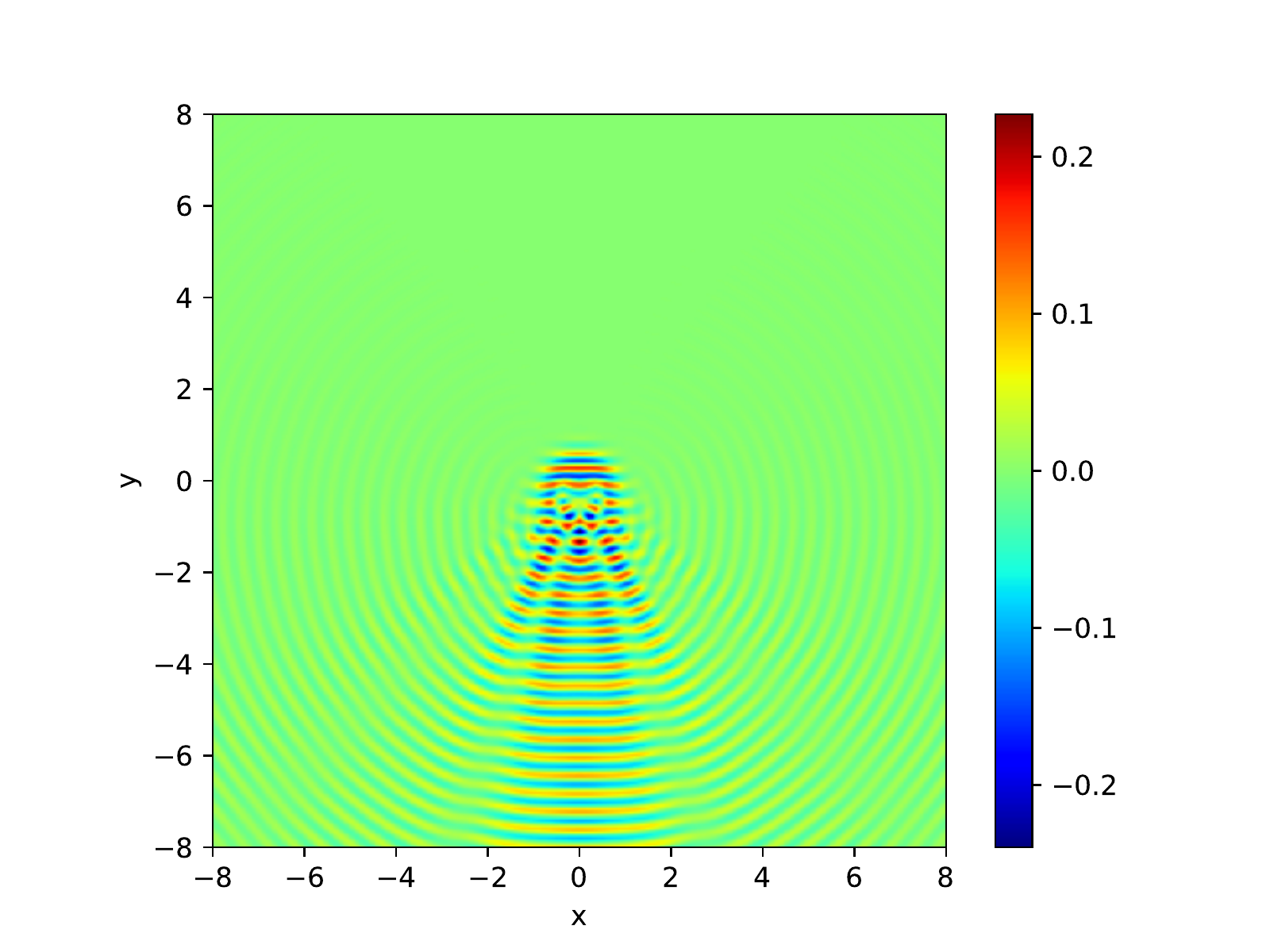}
\end{center}
\end{minipage}

\begin{minipage}[c]{.5\textwidth}
\begin{center}
\includegraphics[width=\textwidth]{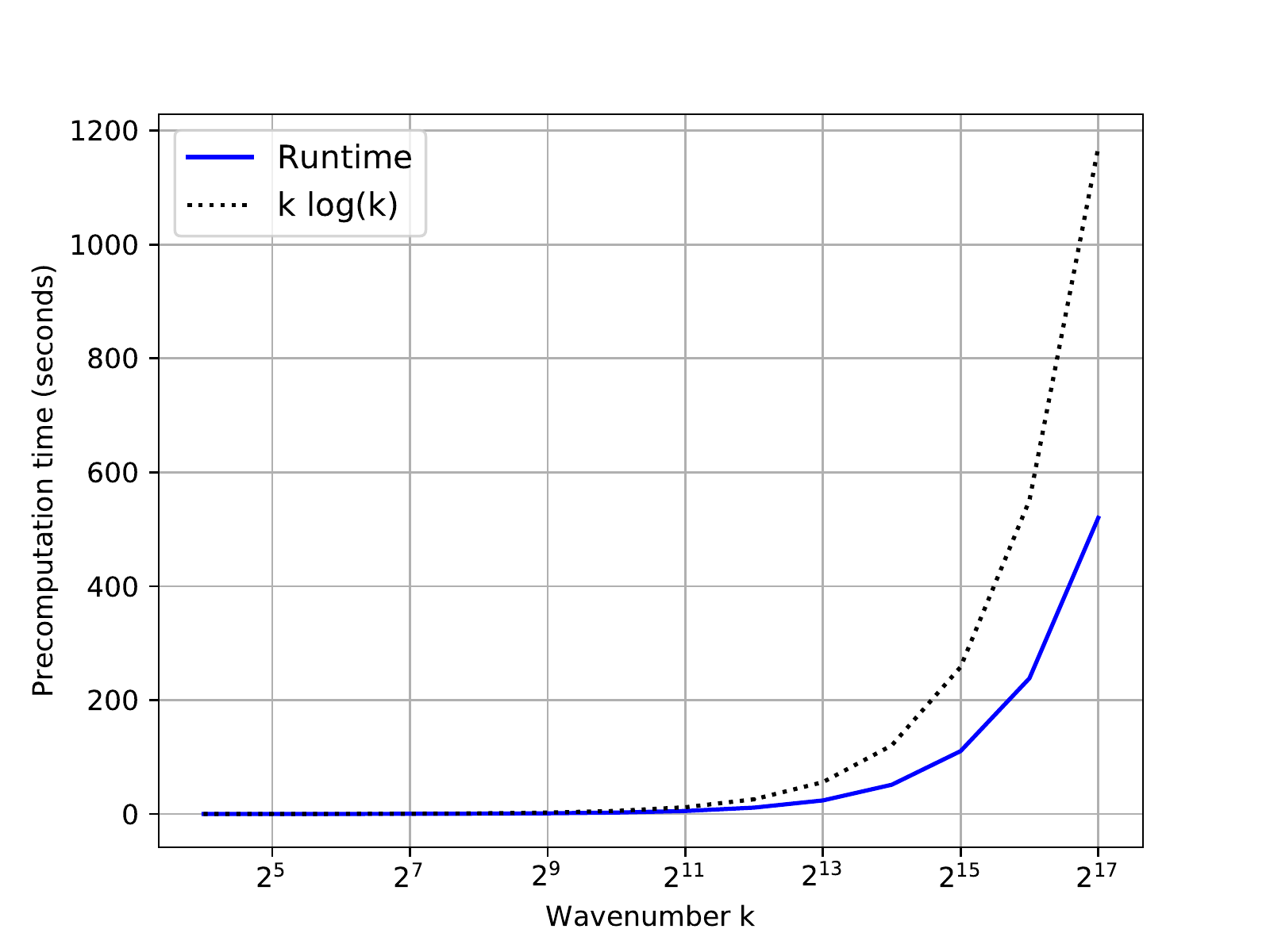}
\end{center}
\end{minipage}
\begin{minipage}[c]{.5\textwidth}
\begin{center}
\includegraphics[width=\textwidth]{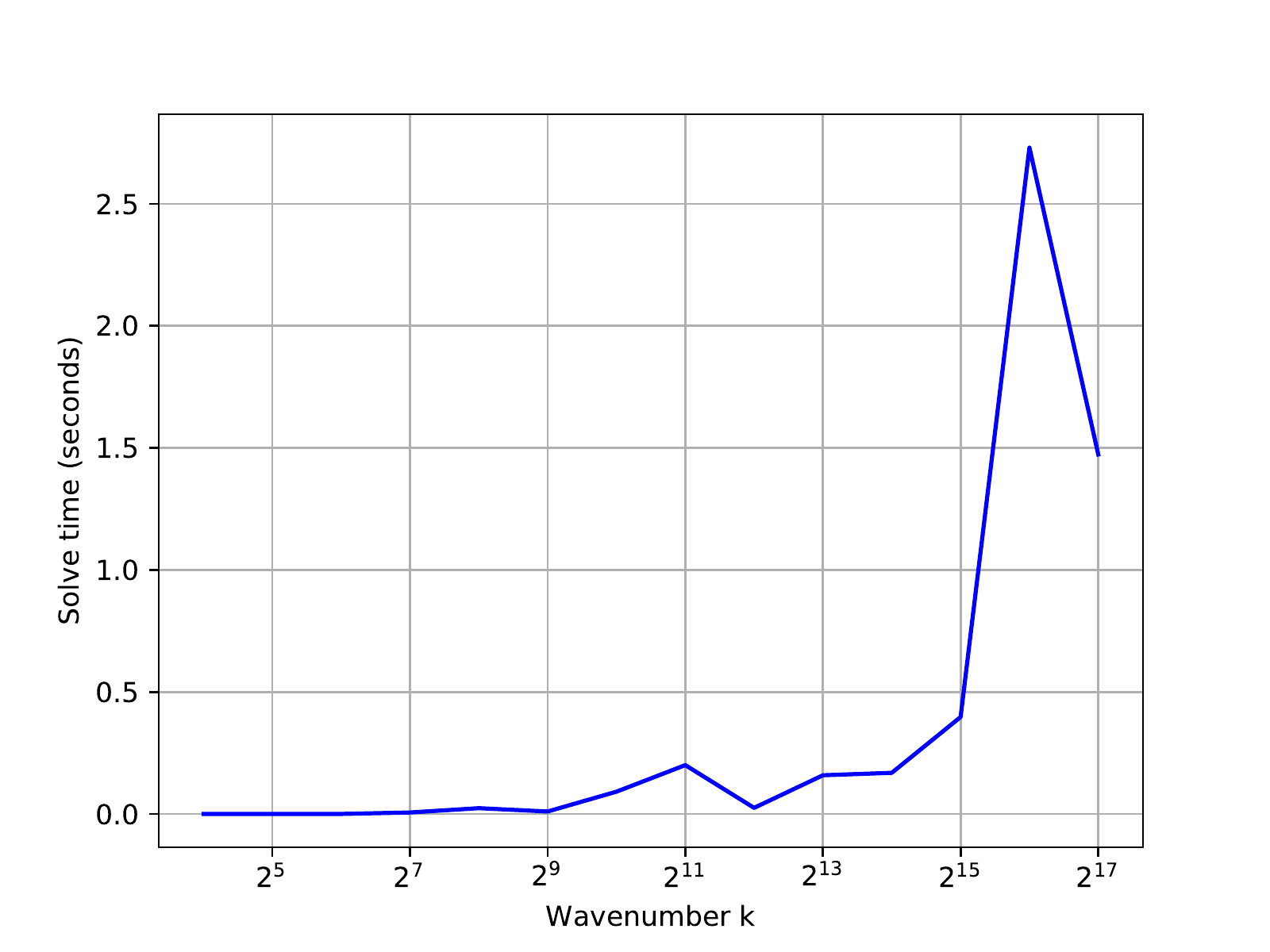}
\end{center}
\end{minipage}

\caption{Figures related to the experiment of Section~\ref{section:experiments:volcano}.
In the upper left is a plot of the function $q(r)$ and in the upper right is
an image of the incident wave when $k=16$.   The image at middle left is of the
real part of the total field $k=16$.  At middle right is an image of
the real part of the scattered field when $k=16$.
At bottom left, the running time of the precomputation phase is plotted as a function of $k$
and at bottom right, the running time of the solution phase is plotted as a function of $k$.
}
\label{figure:volcano1}
\end{figure}

\FloatBarrier

\begin{table}[p]
\begin{center}
\begin{tabular}{crccc}
\toprule
 \addlinespace[.25em]
$k$                                &$m$                                &Maximum absolute                   &Precomp time                       &Solve time                        \\
                                   &                                   &error                              &(in seconds)                       &(in seconds)                       \\
\midrule
 \addlinespace[.25em]
$2^{4}$ &      100 & 9.81\e{-15} & 3.36\e{-02} & 1.74\e{-04}  \\
 \addlinespace[.125em]
$2^{5}$ &      201 & 6.61\e{-14} & 6.36\e{-02} & 2.15\e{-04}  \\
 \addlinespace[.125em]
$2^{6}$ &      402 & 5.96\e{-14} & 1.30\e{-01} & 3.32\e{-04}  \\
 \addlinespace[.125em]
$2^{7}$ &      804 & 2.94\e{-14} & 2.69\e{-01} & 6.22\e{-03}  \\
 \addlinespace[.125em]
$2^{8}$ &     1608 & 4.85\e{-14} & 5.55\e{-01} & 2.38\e{-02}  \\
 \addlinespace[.125em]
$2^{9}$ &     3216 & $\left(6.01\e{-14}\right)$ & 1.16\e{+00}  & 1.02\e{-02}  \\
 \addlinespace[.125em]
$2^{10}$  &     6433 & $\left(9.30\e{-14}\right)$ & 2.47\e{+00}  & 9.14\e{-02}  \\
 \addlinespace[.125em]
$2^{11}$  &    12867 & $\left(1.32\e{-13}\right)$ & 5.24\e{+00}  & 2.00\e{-01}  \\
 \addlinespace[.125em]
$2^{12}$  &    25735 & $\left(5.22\e{-13}\right)$ & 1.11\e{+01}  & 2.51\e{-02}  \\
 \addlinespace[.125em]
$2^{13}$  &    51471 & $\left(1.32\e{-12}\right)$ & 2.38\e{+01}  & 1.58\e{-01}  \\
 \addlinespace[.125em]
$2^{14}$  &   102943 & $\left(3.43\e{-12}\right)$ & 5.13\e{+01}  & 1.68\e{-01}  \\
 \addlinespace[.125em]
$2^{15}$  &   205887 & $\left(1.34\e{-11}\right)$ & 1.10\e{+02}  & 3.97\e{-01}  \\
 \addlinespace[.125em]
$2^{16}$  &   411774 & $\left(5.44\e{-11}\right)$ & 2.38\e{+02}  & 2.73\e{+00}   \\
 \addlinespace[.125em]
$2^{17}$  &   823549 & $\left(1.09\e{-10}\right)$ & 5.20\e{+02}  & 1.47\e{+00}   \\
 \addlinespace[.125em]
\bottomrule
\end{tabular}

\end{center}
\caption{The results of the experiments of Section~\ref{section:experiments:volcano}. Each row
of the table corresponds to one wavenumber $k$ and gives the number $m$ of Fourier modes used
to represent the incident wave, the maximum observed absolute error in the obtained solution
(in cases in which this could be measured), and the time taken by each phase of our solver.
The absolute maximum errors are calculated via comparison with solutions generated
using extended precision arithmetic.  Parentheses are used to indicate cases in which
the wavenumber was too large for the  the extended precision
solution to be verified via a spectral method.
}
\label{table:volcano1}
\end{table}

\begin{figure}[p]
\begin{center}
\includegraphics[width=.75\textwidth]{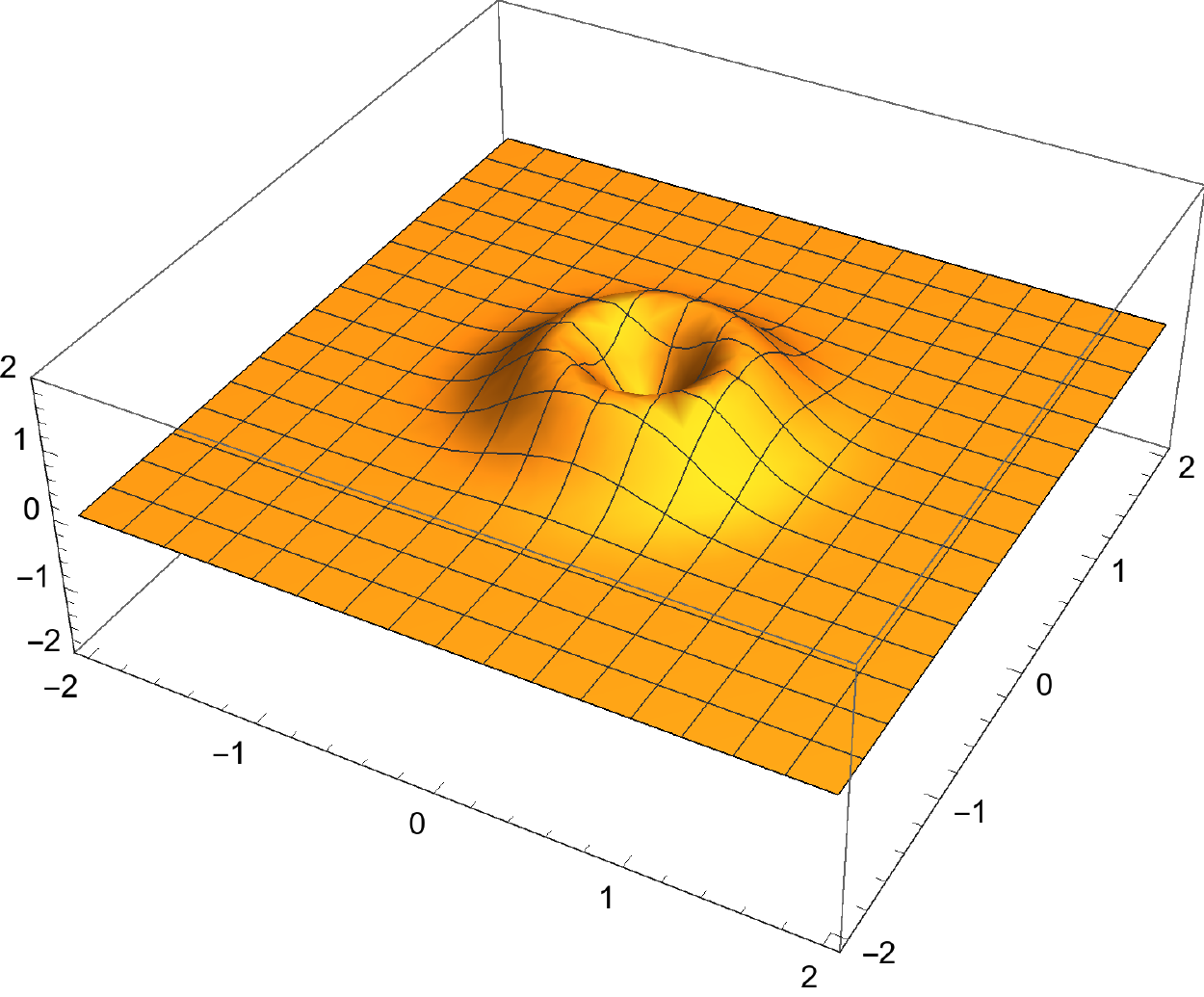}
\end{center}
\caption{A plot of the volcano-shaped 
scattering potential used in the experiment of Section~\ref{section:experiments:volcano}.}
\label{figure:volcano2}
\end{figure}

\FloatBarrier

\begin{figure}[p]
\small
\begin{minipage}[c]{.5\textwidth}
\begin{center}
\includegraphics[width=\textwidth]{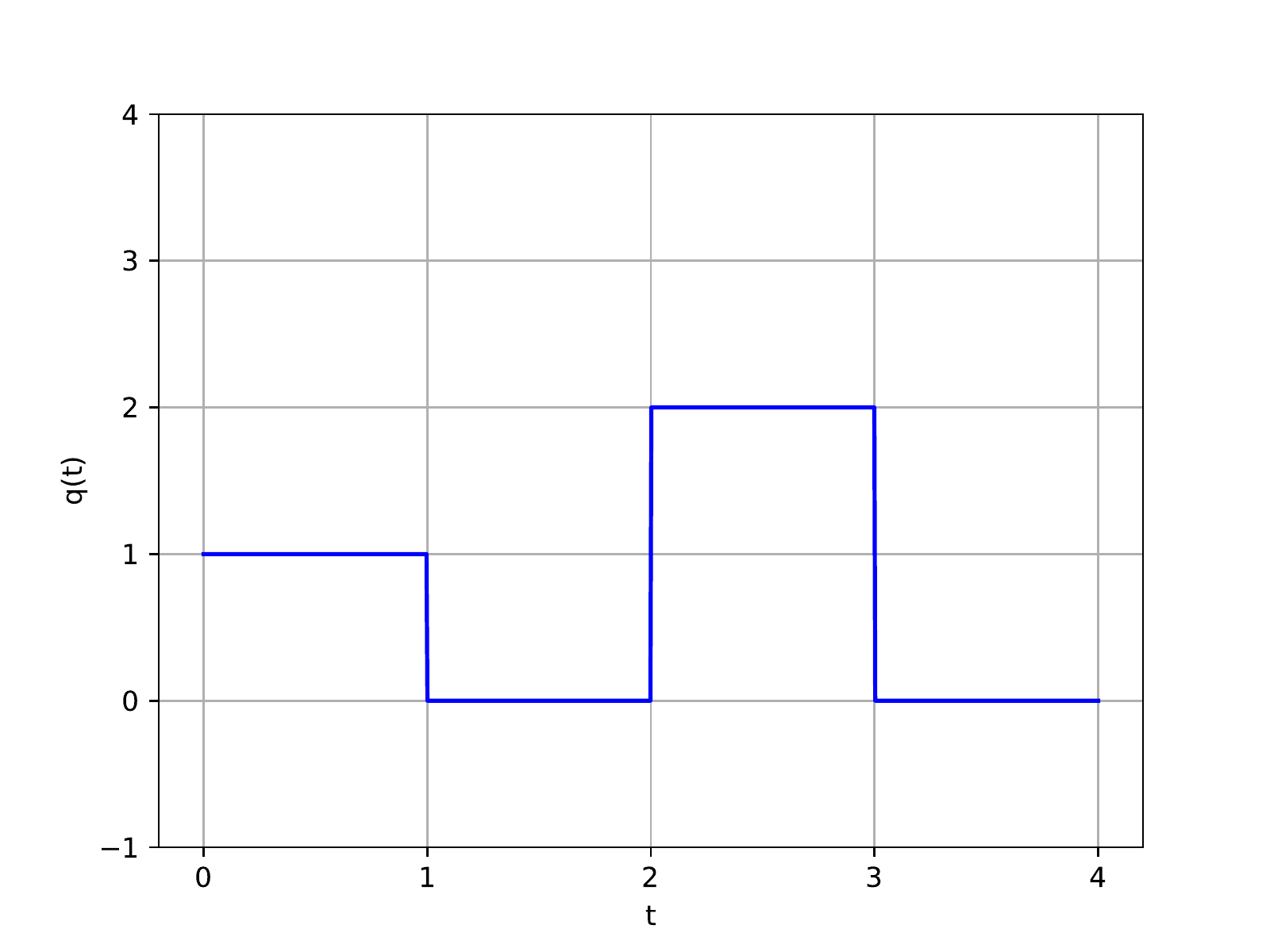}
\end{center}
\end{minipage}
\begin{minipage}[c]{.5\textwidth}
\begin{center}
\includegraphics[width=\textwidth]{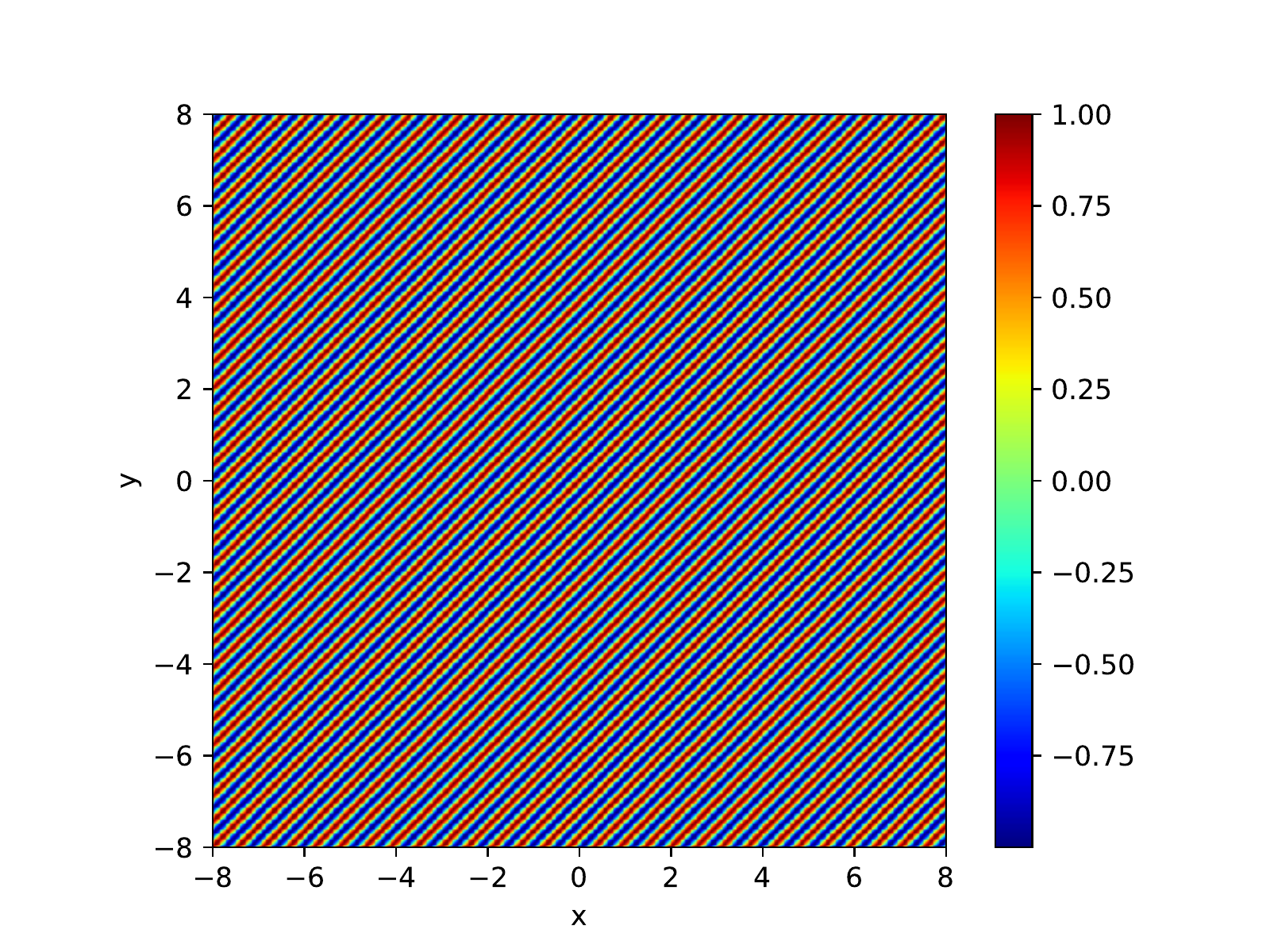}
\end{center}
\end{minipage}

\begin{minipage}[c]{.5\textwidth}
\begin{center}
\includegraphics[width=\textwidth]{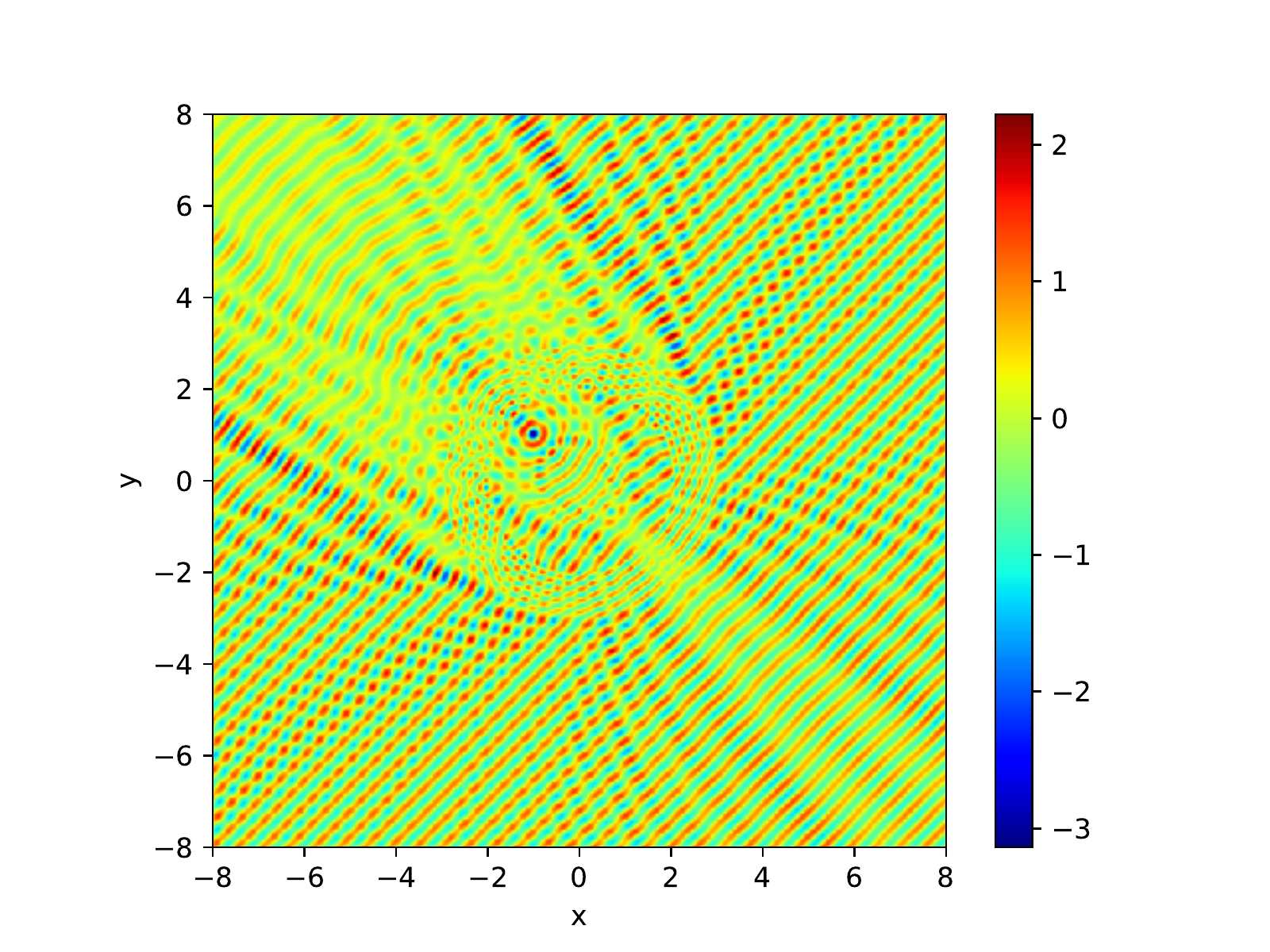}
\end{center}
\end{minipage}
\begin{minipage}[c]{.5\textwidth}
\begin{center}
\includegraphics[width=\textwidth]{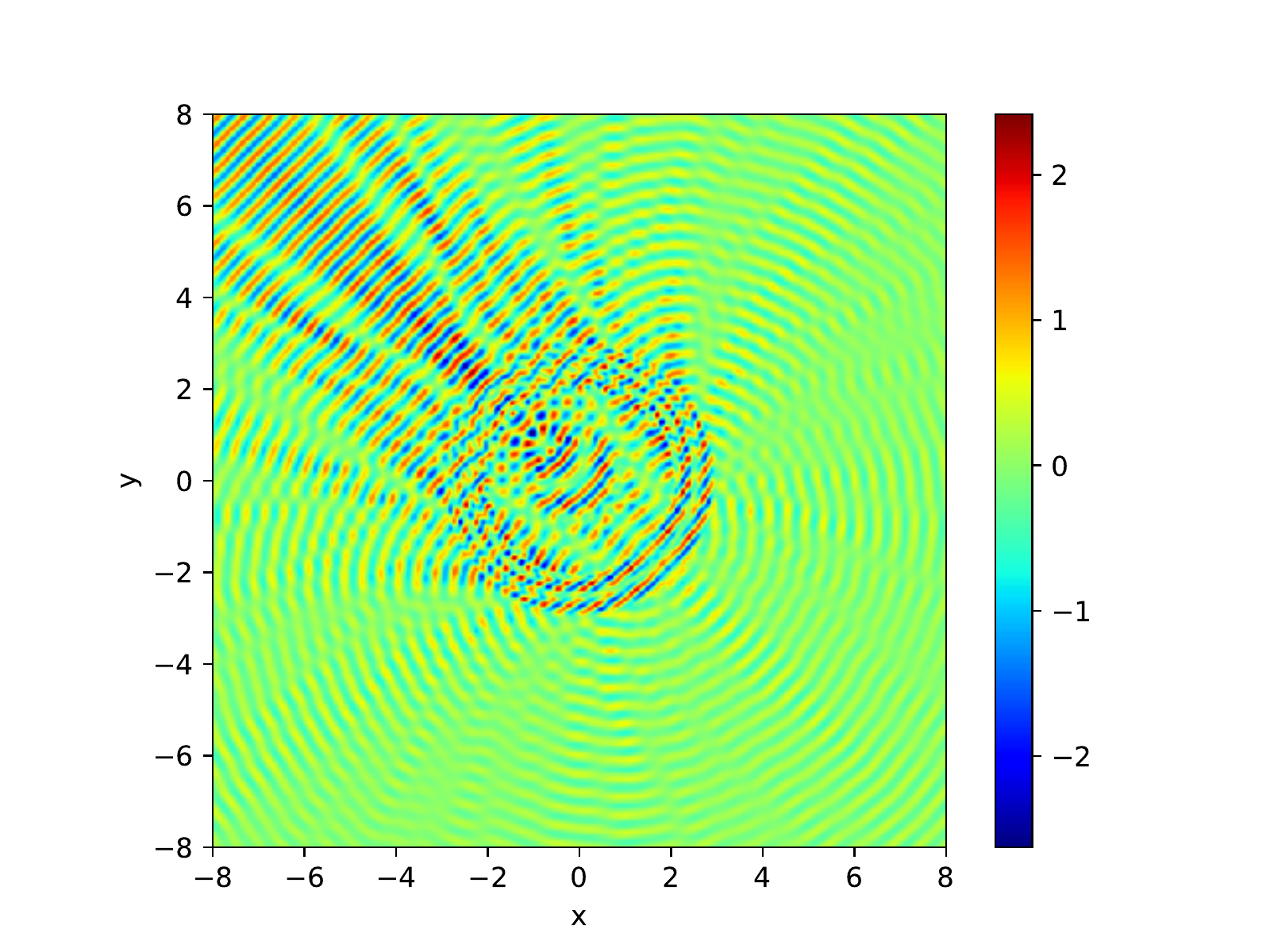}
\end{center}
\end{minipage}

\begin{minipage}[c]{.5\textwidth}
\begin{center}
\includegraphics[width=\textwidth]{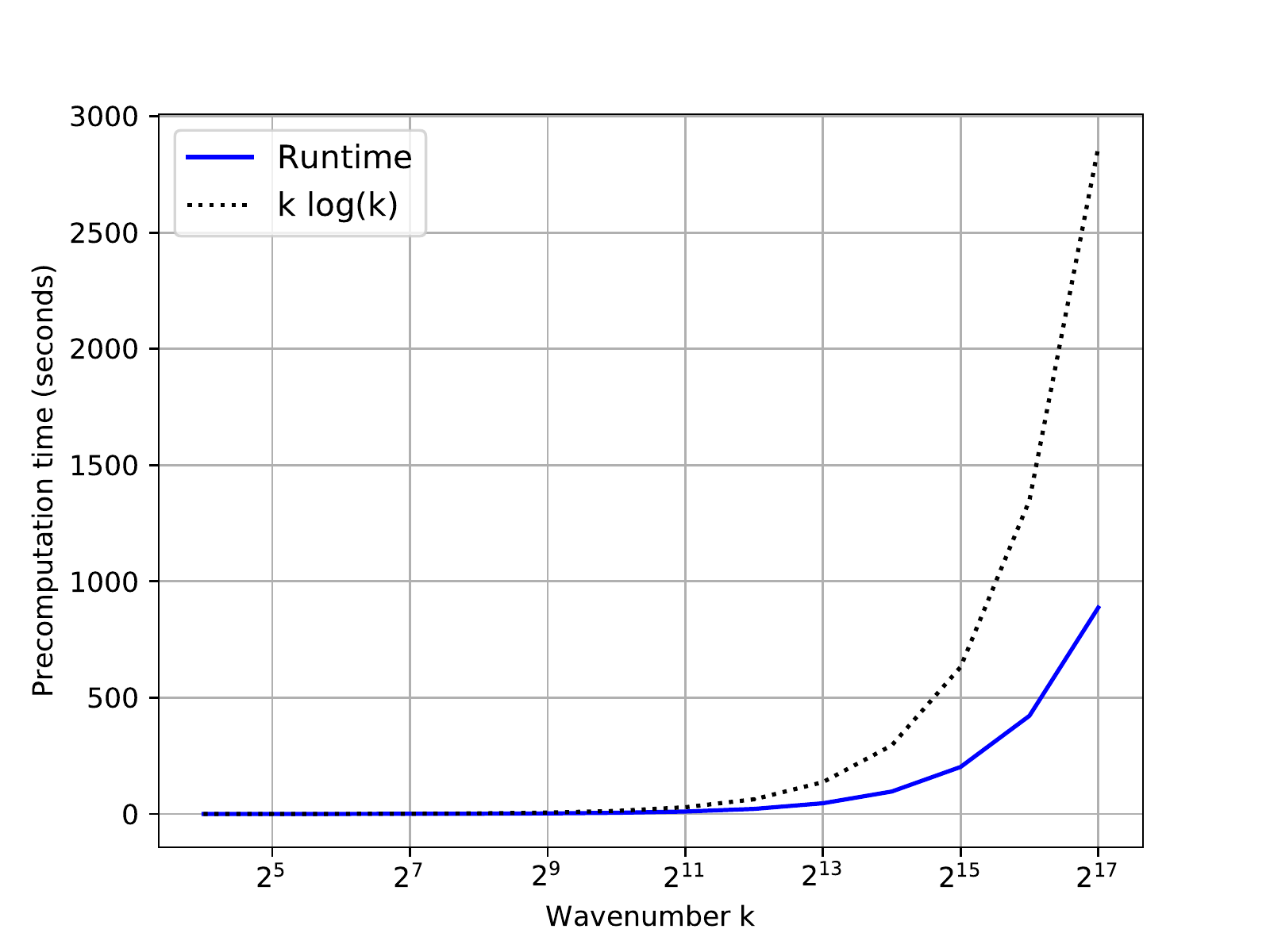}
\end{center}
\end{minipage}
\begin{minipage}[c]{.5\textwidth}
\begin{center}
\includegraphics[width=\textwidth]{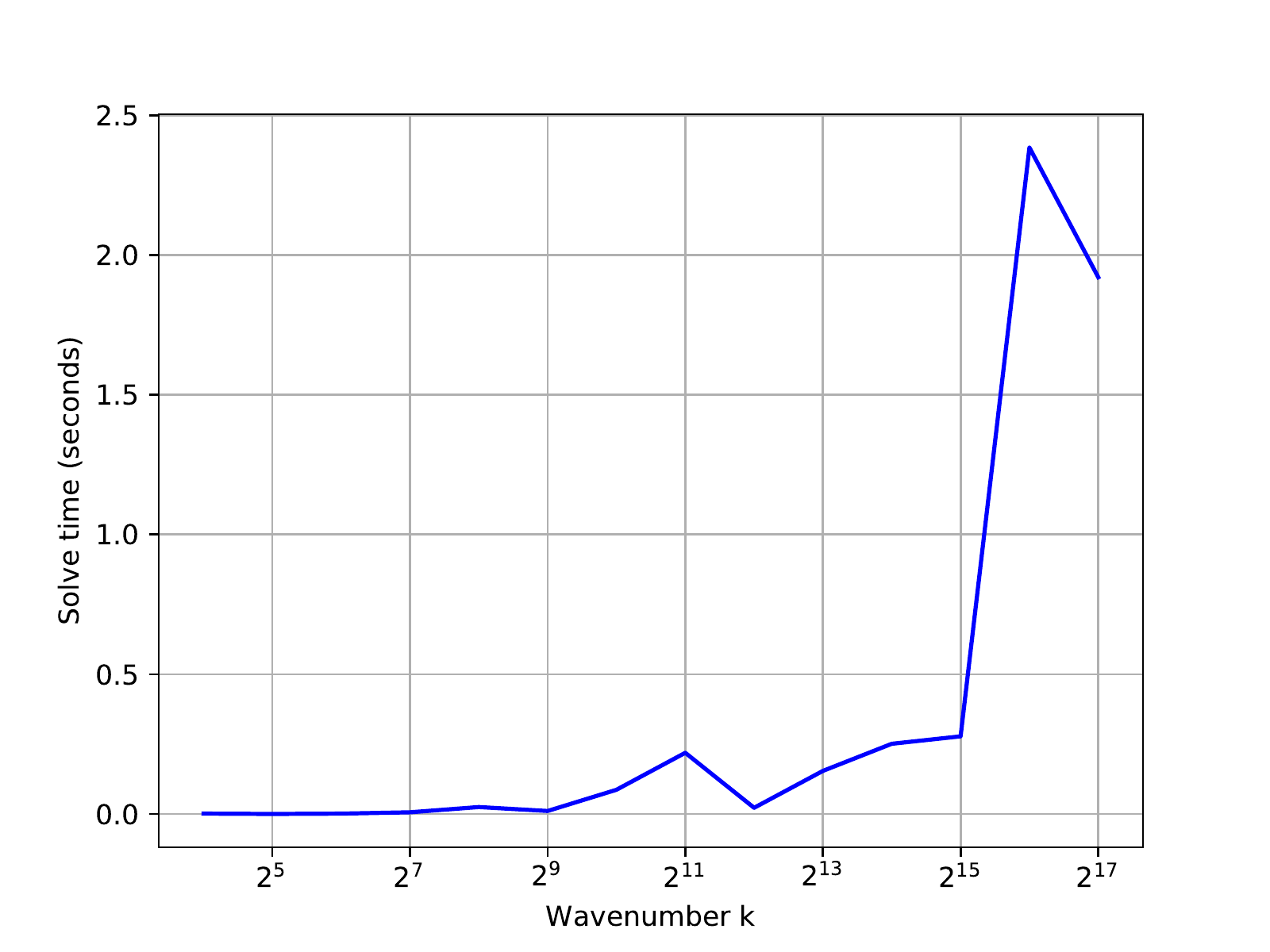}
\end{center}
\end{minipage}

\caption{Figures related to the experiment of Section~\ref{section:experiments:discont}.
In the upper left is a plot of the function $q(r)$ and in the upper right is
an image of the incident wave when $k=16$.   The image at middle left is of the
real part of the total field $k=16$.  At middle right is an image of
the real part of the scattered field when $k=16$.
At bottom left, the running time of the precomputation phase is plotted as a function of $k$
and at bottom right, the running time of the solution phase is plotted as a function of $k$.
}
\label{figure:discont1}
\end{figure}

\begin{table}[p]
\begin{center}
\begin{tabular}{crccc}
\toprule
 \addlinespace[.25em]
$k$                                &$m$                                &Maximum absolute                   &Precomp time                       &Solve time                        \\
                                   &                                   &error                              &(in seconds)                       &(in seconds)                       \\
\midrule
 \addlinespace[.25em]
$2^{4}$ &      100 & 8.03\e{-14} & 8.23\e{-02} & 1.47\e{-03}  \\
 \addlinespace[.125em]
$2^{5}$ &      201 & 1.20\e{-13} & 1.76\e{-01} & 1.85\e{-04}  \\
 \addlinespace[.125em]
$2^{6}$ &      402 & 3.15\e{-13} & 3.81\e{-01} & 1.42\e{-03}  \\
 \addlinespace[.125em]
$2^{7}$ &      804 & 8.20\e{-12} & 7.64\e{-01} & 6.17\e{-03}  \\
 \addlinespace[.125em]
$2^{8}$ &     1608 & 3.46\e{-12} & 1.18\e{+00}  & 2.49\e{-02}  \\
 \addlinespace[.125em]
$2^{9}$ &     3216 & $\left(8.46\e{-12}\right)$  & 2.42\e{+00}  & 1.08\e{-02}  \\
 \addlinespace[.125em]
$2^{10}$  &     6433 & $\left(4.01\e{-11}\right)$  & 5.04\e{+00}  & 8.67\e{-02}  \\
 \addlinespace[.125em]
$2^{11}$  &    12867 & $\left(1.30\e{-10}\right)$  & 1.04\e{+01}  & 2.19\e{-01}  \\
 \addlinespace[.125em]
$2^{12}$  &    25735 & $\left(3.80\e{-10}\right)$  & 2.18\e{+01}  & 2.22\e{-02}  \\
 \addlinespace[.125em]
$2^{13}$  &    51471 & $\left(1.99\e{-09}\right)$  & 4.61\e{+01}  & 1.54\e{-01}  \\
 \addlinespace[.125em]
$2^{14}$  &   102943 & $\left(5.87\e{-09}\right)$  & 9.65\e{+01}  & 2.51\e{-01}  \\
 \addlinespace[.125em]
$2^{15}$  &   205887 & $\left(4.07\e{-08}\right)$  & 2.02\e{+02}  & 2.77\e{-01}  \\
 \addlinespace[.125em]
$2^{16}$  &   411774 & $\left(1.00\e{-07}\right)$  & 4.22\e{+02}  & 2.38\e{+00}   \\
 \addlinespace[.125em]
$2^{17}$  &   823549 & $\left(3.64\e{-07}\right)$  & 8.87\e{+02}  & 1.92\e{+00}   \\
 \addlinespace[.125em]
\bottomrule
\end{tabular}

\end{center}
\caption{The results of the experiments of Section~\ref{section:experiments:discont}. Each row
of the table corresponds to one wavenumber $k$ and gives the number $m$ of Fourier modes used
to represent the incident wave, the maximum observed absolute error in the obtained solution
(in cases in which this could be measured), and the time taken by each phase of our solver.
The absolute maximum errors are calculated via comparison with solutions generated
using extended precision arithmetic.  Parentheses are used to indicate cases in which
the wavenumber was too large for the  the extended precision
solution to be verified via a spectral method.
}
\label{table:discont1}
\end{table}

\begin{figure}[p]
\begin{center}
\includegraphics[width=.75\textwidth]{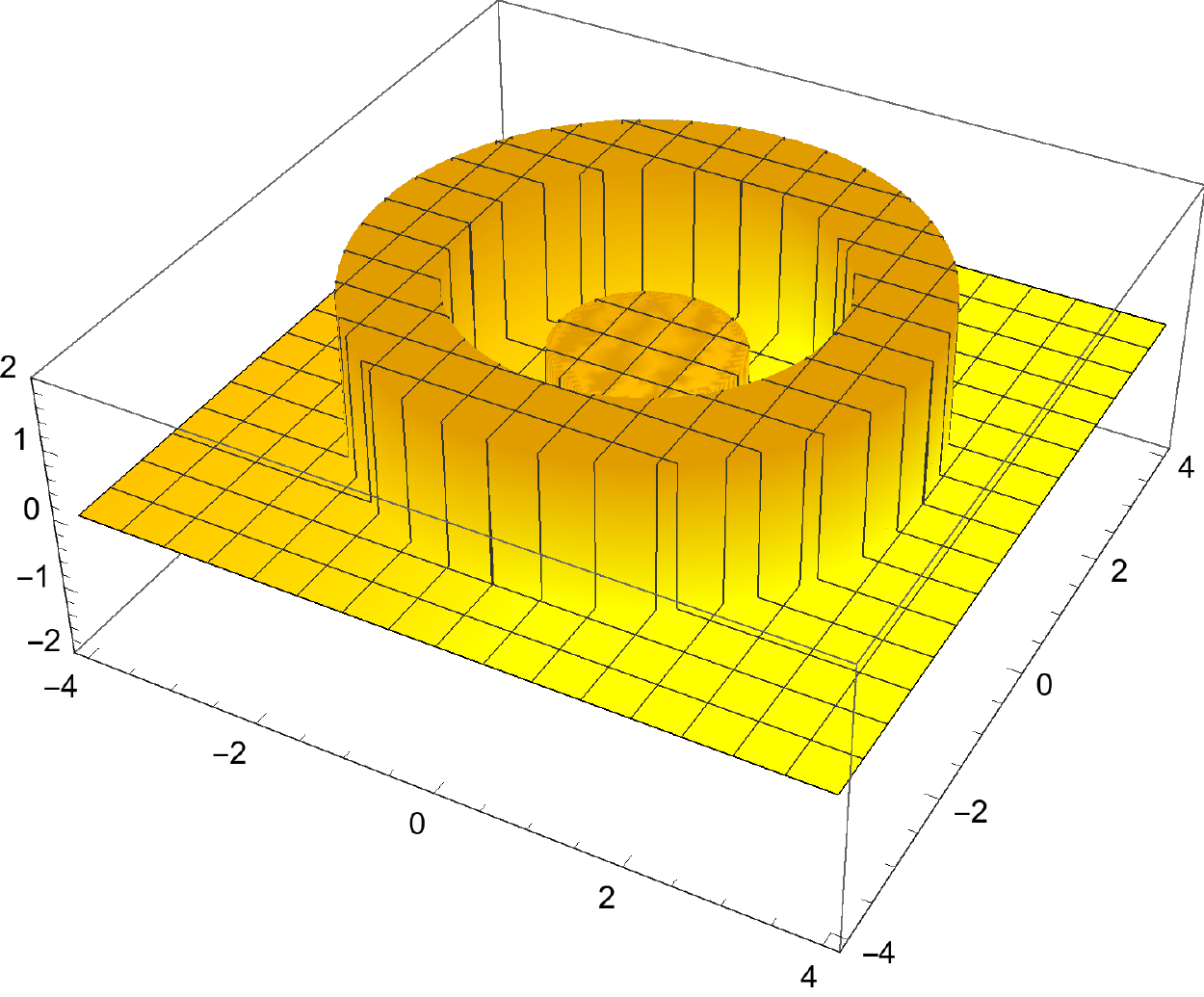}
\end{center}
\caption{A plot of the discontinuous 
scattering potential used in the experiment of Section~\ref{section:experiments:discont}.}
\label{figure:discont2}
\end{figure}





  



\end{document}